\documentclass[a4paper,11pt,reqno,noindent]{amsart}
\usepackage[centertags]{amsmath}
\usepackage{amsfonts,amssymb,amsthm,amscd,dsfont,cases,esint,enumerate,color}
\usepackage{comment}
\usepackage[T1]{fontenc}
\usepackage[english]{babel}
\usepackage[applemac]{inputenc}
\usepackage[body={15cm,21.5cm},centering]{geometry} 
\usepackage{fancyhdr}
\numberwithin{equation}{section}

\usepackage{graphicx}
\usepackage{esint}
\usepackage{amsmath}
\usepackage{amsfonts}
\usepackage{amssymb}
\usepackage{amsthm}
\usepackage{comment}
\usepackage{appendix}
\usepackage{enumitem}
\usepackage{mathtools}

\theoremstyle{plain}
\newtheorem{theor10}{Theorem}

\newtheorem{prop10}[theor10]{Proposition}

\newtheorem{cor10}[theor10]{Corollary}

\newtheorem{lem10}[theor10]{Lemma}

\newtheorem{theor0}{Theorem}[section]

\newtheorem{lem0}[theor0]{Lemma}
\newenvironment{lem}
  {\pushQED{\qed}\begin{lem0}}
  {\popQED\end{lem0}}
\newtheorem{prop0}[theor0]{Proposition}

\newtheorem{cor0}[theor0]{Corollary}

\theoremstyle{definition}
\newtheorem{rems0}[theor0]{Remarks}

\newtheorem{rem0}[theor0]{Remark}

 \newtheorem{hypo0}[theor0]{Hypothesis}

\theoremstyle{plain}
\newtheorem{as0}[theor0]{Assumption}

\newtheorem*{asn0*}{\assumptionnumber}
  \providecommand{\assumptionnumber}{}
\makeatletter
\newenvironment{asn0}[2]
   {\renewcommand{\assumptionnumber}{Assumption \!#1 {\normalfont--- #2}}
    \begin{asn0*}
    \protected@edef\@currentlabel{{\normalfont#1}}}
   {\end{asn0*}}
\makeatother

\makeatletter
\newenvironment{asn01}[1]
   {\renewcommand{\assumptionnumber}{Assumption \!#1}
    \begin{asn0*}
    \protected@edef\@currentlabel{{\normalfont#1}}}
   {\end{asn0*}}
\makeatother

\newcommand{\N}{\mathbb N}

\newcommand{\R}{\mathbb R}
\newcommand{\Z}{\mathbb Z}

\newcommand{\Id}{\operatorname{Id}}

\newcommand{\vp}{\varphi}

\newcommand{\ZZ}{\mathbb{Z}}
\newcommand{\BB}{\mathcal{B}}

\newcommand{\calz}{\mathbb{B}}
\newcommand{\RR}{\mathbb{R}}

\newcommand{\spd}{\spadesuit}

\newcommand{\maxi}{\mathcal{M}}
\newcommand{\rd}{{r_\diamond}}
\newcommand{\di}{{\diamond}}

\newcommand{\expecm}[1]{\mathbb{E}\big[ #1 \big]}

\newcommand{\expecb}[1]{\mathbb{E}\Big[#1\Big]}

\newcommand{\expec}[1]{\mathbb{E}[#1]}

\newcommand{\Cc}{\mathcal C}

\newtheorem{theorem}{Theorem}
\newtheorem{proposition}{Proposition}
\newtheorem{corollary}{Corollary}
\newtheorem{lemma}{Lemma}

\usepackage[colorlinks,citecolor=black,urlcolor=black]{hyperref}

\title[Corrector bounds for the degenerate random conductance model]{Moment bounds on correctors for the degenerate random conductance model}

\author[A. Gloria]{Antoine Gloria}
\address[Antoine Gloria]{Sorbonne Universit\'e, CNRS, Universit\'e de Paris, Laboratoire Jacques-Louis Lions, 75005~Paris, France  \& Universit\'e Libre de Bruxelles, D\'epartement de Math\'ematique, 1050~Brussels, Belgium}
\email{antoine.gloria@sorbonne-universite.fr}

\author[S. Qi]{Siguang Qi}
\address[Siguang Qi]{Sorbonne Universit\'e, CNRS, Universit\'e de Paris, Laboratoire Jacques-Louis Lions, 75005~Paris, France }
\email{siguang.qi@sorbonne-universite.fr}

\begin{document}
\selectlanguage{english}

\maketitle

\begin{abstract}
We study the random conductance model on the lattice $\Z^d$, i.e. we consider a linear,
finite-difference, divergence-form operator with random conductances $a$. We allow the conductances $a$ to be unbounded and degenerate. Assuming the conductances satisfy a spectral-gap inequality, we establish sharp bounds on the spatial growth of correctors, together with a quantitative relation between the stochastic integrability of the correctors and that of $a$. 

\bigskip\noindent
{\sc MSC-class:} 35R60, 35B27, 35B65, 60H07

\end{abstract}

\setcounter{tocdepth}{1}
\tableofcontents

\section{Introduction}

This article is devoted to establishing moment bounds on correctors in the degenerate random conductance model. In the first part of this introduction, we discuss moment bounds on correctors in stochastic homogenization in general (and more specifically, in the continuum setting). In the second part, we state our main results for degenerate conductances on the lattice and compare them with the existing literature.

\subsection{Stochastic integrability in stochastic homogenization}

Let $A$ be a symmetric non-negative stationary random matrix field on $\R^d$. Homogenization of linear elliptic equations in divergence form amounts to understanding the large-scale behavior of the operator $-\nabla \cdot A \nabla$. A key quantity in this theory is the corrector, that is, the correction to a linear function that makes it $A$-harmonic. This is natural, because smooth functions are locally linear, so solutions of $-\nabla \cdot A \nabla u=f$ must locally look like a linear function plus a corrector. In probabilistic approaches based on the environment viewed by the particle, the corrector is precisely the object that turns the process into a martingale. The accuracy of this approximation (the so-called two-scale expansion) is measured by two parameters: a convergence rate, which is mostly related to the correlations of $A$ (and how they are measured), and a stochastic-integrability estimate, which captures finer concentration properties. The strongest possible versions of these two parameters are the central limit theorem (CLT) scaling and Gaussian stochastic integrability. Another feature of $A$, however, plays a decisive role for the latter: the possible degeneracy of the coefficient field.

\medskip

Let us be more precise. Given such a field $A$ and a unit direction $e\in \R^d$, we define the corrector $\phi$ in direction $e$ as a solution in $\R^d$, unique up to an additive constant, of 
\[
-\nabla \cdot A (\nabla \phi+e)=0,
\] 
and assume for the moment that $\nabla \phi$ is a centered stationary field with finite first moment.
By the ergodic theorem, $\fint_{B_R} \nabla \phi$ converges to $\expecm{\nabla \phi}=0$ as $R \uparrow +\infty$ almost surely. Assuming in addition that $A$ is weakly correlated, we expect $\fint_{B_R} \nabla \phi$ to decay at the CLT scaling, so that the sequence of non-negative numbers
\[
\Cc_R:=R^{d/2} \Big| \fint_{B_R} \nabla \phi\Big|
\]
is expected to be controlled uniformly in $R$ (here $B_R$ denotes the ball of radius $R$ centered at the origin). An important question is therefore the stochastic integrability of $\Cc_R$, and more specifically how it depends on the stochastic integrability of $A$ in degenerate models.
The following table summarizes the state of the art on this question together with the results of the present contribution, which we expect to be optimal. For simplicity, we take $A=a \Id$ in this table:

\medskip
\begin{center}
\begin{tabular}{|c|c|c|c|}
\hline
 Stochastic   & Mixing  & Stochastic   & Reference
\\
 integrability  of $a$ & condition & integrability of $\Cc_R$ & 
\\
\hline
 $0<\lambda \le a \le 1$ & FRD  & $\expecm{\exp(\frac1{C_\lambda} \frac {\Cc_R^2}{\log^{\alpha(d)}(2+\Cc_R)})}\le 2$ & \cite{GO-25} 
\\
\hline
 $\expecm{\exp(\frac1C \log^{\frac2\alpha} (a+a^{-1}))} \le 2$
& LSI  & $\expecm{\exp(\frac1C \log^{\frac2\alpha} \Cc_R )}\le 2$ & \cite{CGQ-24}
\\
\hline
 $\expec{a^{p}+a^{-p}}<\infty$ & i.i.d. &  $\sup_{R\ge 1} \expec{\Cc_R^{\frac1C p}}<\infty$ & Theorem~\ref{thm.main}.
\\
\hline
\end{tabular}
\end{center}

\medskip

\begin{itemize}
\item The first model is a uniformly elliptic and bounded coefficient field $a$, that is, $0<\lambda\le a\le 1$, with finite range of dependence (FRD). In this case, the stochastic integrability of $\Cc_R$ is Gaussian up to a poly-logarithmic correction $0<\alpha(d)<\infty$; see~\cite{GO-25} by the first author and Otto (see also \cite{GO4}, and  \cite{AKM2} by Armstrong, Kuusi, and Mourrat for strictly sub-Gaussian stochastic integrability).
\item The second model is a log-normal coefficient of the form $a(x)=\exp\Big(G(x)(1+|G(x)|^2)^{\frac{\alpha-1}{2}} \Big)$ with $0<\alpha<2$ and $G$ a stationary centered Gaussian field with integrable covariance.  Such a model satisfies a logarithmic Sobolev inequality (LSI), and the stochastic integrability of $\Cc_R$ is of the same type as that of $a$ (with a different constant $C$); see~\cite{CGQ-24} by Clozeau and the present authors.
\item The model treated in the present contribution is even more degenerate: we only assume finite moments of $a$ and $a^{-1}$. We prove that there exists a constant $C \ge 1$ such that, for every sufficiently large $p$, if the former have finite $p$-th moments, then $\Cc_R$ has finite moment of order $\frac1C p$; see Theorem~\ref{thm.main} below. Motivated by the probabilistic literature, we work with the discrete version on $\Z^d$, namely the degenerate random conductance model with independent and identically distributed (i.i.d.) coefficients. We also refer to  \cite{MR3949105} by Andres and Neukamm for quantitative (yet, less precise) results on this model.
\end{itemize}

\medskip

At this stage, we should also mention another popular degenerate model: the discrete Laplacian on the percolation cluster (in the supercritical regime). Armstrong and Dario \cite{AD-18}, and later Dario  \cite{Dario},
showed that $\Cc_R$ has stretched-exponential stochastic integrability (optimality is unclear).

\medskip

Another important question, which is to some extent a prerequisite for the present work, is under which assumptions on the conductances correctors exist and qualitative homogenization holds. This has been a very active field of probability theory, in particular, in the study of random walks in random environments. Our assumptions on random conductances will be stronger than the minimal assumptions for qualitative homogenization. In particular, for stationary and ergodic random conductances, the state of the art is 
\begin{equation}\label{e.best-pq}
\expec{a^q}+\expec{a^{-p}} <\infty, \quad \frac1q+\frac1p <\frac2{d-1}, \quad p,q>1,
\end{equation}
which was recently proved to be sufficient for the invariance principle by Bella and Sch\"affner in \cite{MR4079437}, building and improving on the result by Andres, Deuschel, and Slowik in \cite{AndresDeuschelSlowik2015}.
When the law of the conductances is i.i.d., the results are stronger, and it suffices that the set of conductances with value larger than a fixed positive constant percolates, cf.~\cite{MR3078279} by Andres, Barlow,   Deuschel, and Hambly.
We shall place ourselves in the setting when correctors are well-defined, and shall rely on~\cite{bella2018liouville} by Bella, Fehrman, and Otto.

\subsection{Assumptions and main results}\label{sec:main-results}

In this work we consider the random conductance model, that is, let $(a_e)_{e\in\calz^d}$ be a family of i.i.d. random variables defined on $\calz^d$, the set of all edges in $\ZZ^d$, satisfying the following moment bound
\begin{equation}\label{eq.mombd-a}
\expec{a^\gamma}^\frac1\gamma+\expec{a^{-\gamma}}^\frac1\gamma\eqqcolon\Gamma<+\infty
\end{equation}
for some\footnote{We only need $\gamma>d$ in the argument, but since the result requires $\gamma$ to be large enough anyway, we keep the convenient choice $\gamma\ge d+1$ in what follows.} $\gamma \ge d+1$ (that is, in the setting of \eqref{e.best-pq}, ``$p=q=\gamma$, and $\frac1p+\frac1q=\frac2\gamma\le \frac2{d+1} <\frac2{d-1}$''), where $\expec{\cdot}$ denotes the expectation. For each unit direction $e_i$ in $\ZZ^d$, we consider the following corrector equation 
\begin{equation}\label{eq.corrector}
-\nabla^*\cdot A (\nabla\phi_i+e_i)=0,
\end{equation}
where $A:\ZZ^d\to \RR^{d\times d}$ is a matrix defined at vertices $x$ via
\[
A(x)=\operatorname{diag}(a_{\{x,x+e_1\}},a_{\{x,x+e_2\}},\cdots,a_{\{x,x+e_d\}}),
\]
$\nabla$ is the forward discrete gradient, and $\nabla^*$ the backward gradient -- see precise notation below. 
 The flux corrector $\sigma_i$ in direction $e_i$ is a skew-symmetric matrix field whose entry $\sigma_{ijk}$ satisfies on $\Z^d$ 
\begin{equation}\label{eq.flux}
-\Delta \sigma_{ijk}= \nabla_{j} \big(A (\nabla\phi_{i}+{e_i})\big)_k-\nabla_{k} \big(A (\nabla\phi_{i}+{e_i})\big)_j.
\end{equation}
These objects are well-defined in the following sense, see \cite{bella2018liouville,AndresDeuschelSlowik2015}.
\begin{lem}\label{lem:liouville}
Let  $a$ satisfy \eqref{eq.mombd-a} for some $\gamma\ge d+1$.
There exist a vector field $\phi=\{\phi_i\}_{1\le i \le d}$ and a tensor field $\sigma=\{\sigma_{ijk}\}_{1\le i,j,k\le d}$ with the following properties. For each $1\le i,j,k\le d$, $\phi_i$ solves \eqref{eq.corrector} and $\sigma_{ijk}$ solves \eqref{eq.flux}. The gradient fields are stationary (that is, their laws are invariant under translations), have some finite moments, and have zero mean: $\expec{\nabla \phi_i}=\expec{\nabla \sigma_{ijk}}=0$, and 
\[
\sum_{i=1}^d \expecb{a |\nabla \phi_i|^2}^\frac12+\sum_{i=1}^d \expecb{|\nabla \phi_i|^\frac{2\gamma}{\gamma+1}}^\frac{\gamma+1}{2\gamma} + \sum_{i,j,k=1}^d \expecb{|\nabla \sigma_{ijk}|^\frac{2\gamma}{\gamma+1}}^\frac{\gamma+1}{2\gamma} \,\lesssim \, \Gamma.
\]
Moreover, the correctors grow sub-linearly at infinity:
\begin{equation}\label{est.sublinear}
\lim_{n\to+\infty} \frac1n \fint_{B_{n}} |(\phi_i,\sigma_{ijk})|=0,
\end{equation} 
where $\fint_{B_{n}}$ denotes the average on the ball $B_n$ of size $n$ (see notation for details).
\end{lem}
Our main result is as follows. 
\begin{theorem} \label{thm.main}
Let  $a$ satisfy \eqref{eq.mombd-a} for some $\gamma\ge d+1$.
There exist constants $C_0, C_1, C_2 \ge 1$, depending on dimension and on the effective elliptic ratio $\Lambda\coloneqq\expec{a^{d+1}}^{\frac1{ d+1 }}\expec{a^{-(d+1)}}^{\frac1{ d+1 }}\ge 1$ only, such that if  $\gamma \ge C_2$, then
\begin{equation*}
\sup_{R\ge 1} \expecb{\Big|R^{\frac d2} \fint_{B_R}\nabla (\phi,\sigma) \Big|^{\frac\gamma {C_1}}}^{\frac {C_1}{\gamma}}\lesssim_{\Lambda} (\gamma\Gamma)^{C_0},
\end{equation*}
where $\Gamma$ is defined in \eqref{eq.mombd-a}.
\end{theorem}
Theorem~\ref{thm.main} shows that the loss between the stochastic integrability $\gamma$ of $a$ and that $\frac\gamma{C_1}$ of the corrector is at most linear. This is new and sharp, as the one-dimensional case already shows: By an explicit computation, $C_1=C_2=1$ for $d=1$. In dimensions $d>1$, we only prove the existence of finite constants $C_1,C_2\ge 1$; we do not attempt to optimize them here. 
All in all, the degeneracy appears twice: in the left-hand side through $\gamma$ (correctors cannot be more integrable than the coefficients), and  in the right-hand side through $\Gamma$. The latter dependence with respect to $\Gamma$ is also quantitative: this is crucial in order to use it for more integrable coefficients (like log-normal fields). The constant $C_0$ depends on the strategy of proof, and it cannot be sharp (there is necessarily a loss when using spectral-gap -- see \cite{GNO-reg}).
As a direct corollary of the above theorem, we obtain bounds on correctors. 
\begin{corollary}\label{cor.stationary}
Under the assumptions of Theorem~\ref{thm.main},  if $\gamma \ge C_2$, the correctors satisfy, for all $x\in \Z^d$,
    \begin{equation}\label{eq.growth-corrector}
        \expecb{|(\phi,\sigma)(x)-(\phi,\sigma)(0)|^{\frac\gamma{C_1}}}^\frac{C_1}{\gamma}\lesssim_{ \Lambda }( \gamma\Gamma)^{C_0} \left\{
        \begin{aligned}
            &\sqrt{ |x|},\quad &\text{ if }d=1;\\
            &\log^\frac12(1+|x|),\quad &\text{ if }d=2;\\
            &1,\quad &\text{ if }d\ge3.
        \end{aligned}\right.
    \end{equation}
Hence, if $d\ge3$, there is a unique stationary solution $(\phi,\sigma)$ to the equations \eqref{eq.corrector} and \eqref{eq.flux} such that $\expec{(\phi,\sigma)}=0$. In this case, 
    \begin{equation}\label{eq.correctorbound}
        \expec{|(\phi,\sigma)|^{\frac\gamma{C_1}}}^{\frac{C_1}\gamma}\lesssim_{ \Lambda }( \gamma\Gamma)^{C_0}.
    \end{equation}
\end{corollary}

\subsection{Comparison to the literature and general strategy of the proof}

In \cite{MR3949105}, Andres and Neukamm proved the first quantitative bounds on correctors for the degenerate random conductance model, based on the semigroup method of \cite{GNO1}. Their result is restricted to dimension $d\ge3$, the results on the flux corrector $\sigma$ are not sharp, and the relation between the integrability of the correctors and the coefficients is not explicit -- which is our main purpose in this work. The quantitative forms $\frac\gamma{C_1}$ and $(\gamma\Gamma)^{C_0}$ in the estimate of Theorem~\ref{thm.main} are indeed robust enough to cover less degenerate regimes as well. In particular, if $A$ is uniformly elliptic or has stretched-exponential moments, the theorem implies stretched-exponential integrability for $\Cc_R$ (with a non-optimal exponent -- which significantly improves the recent work \cite{BellaKniely} by Bella and Kniely). In the log-normal case, it yields the sharp log-normal integrability (recovering \cite{CGQ-24}). 

In the work \cite{AndresHalberstam2021} by Andres and Halberstam, algebraic integrability assumptions on $a$, together with a variant of the spectral-gap inequality, are used to obtain quantitative bounds on the random scale above which the heat kernel associated with $-\nabla^*\!\cdot A\nabla$ satisfies lower Gaussian estimates of the same form as in the uniformly elliptic case. In particular, obtaining such quantitative information requires stronger stochastic integrability of the coefficient field than is needed for the qualitative invariance principle. A similar loss appears in our work through the constant $C_1$ in Theorem~\ref{thm.main}. 

Before we describe the general strategy of our proof, let us mention the recent work \cite{AK-25} by Armstrong and Kuusi, which introduces a concept of large-scale ellipticity and a general method to deal with degenerate coefficients. It is not clear to us whether this approach gives some insight into the values of the best constants $C_0,C_1,C_2$, which remain the main open issue in Theorem~\ref{thm.main}.

\medskip

One aim of the present contribution, besides the results themselves, is to show that the elementary method introduced and developed by Clozeau and the present authors in \cite{CGQ-24} can be pushed to treat the degenerate random conductance model.

\medskip

More precisely, our argument relies on two ingredients:
\begin{enumerate}
\item a moment assumption $\expec{a^{\gamma}+a^{-\gamma}}<\infty$ for some $\gamma\ge d+1$;
\item a spectral-gap inequality for the law of $a$.
\end{enumerate}
Starting from these assumptions, the proof proceeds in five steps:
\begin{enumerate}
\item define an effective length scale $r_\diamond\ge2$ such that, for all $r\ge r_\diamond$, the averages of $a^{d+1}$ and $a^{-(d+1)}$ over $B_r$ are bounded by deterministic multiples of their expectations;
\item derive large-scale quenched Meyers estimates at scale $r_\diamond$ for $-\nabla^*\cdot A\nabla$;
\item convert these estimates into large-scale hole-filling and annealed Meyers estimates, whose stochastic integrability is inherited from that of $r_\diamond$;
\item define a minimal radius $r_\spadesuit\ge2$ such that one has for all $r\ge r_\spadesuit$ that $\fint_{B_r} a|\nabla \phi|^2\le C \fint_{B_{2r}}a$, where $C >\expec{a|\nabla \phi|^2}\expec{a}^{-1}$;
\item combine sensitivity calculus, the spectral-gap inequality, and Caccioppoli-type estimates to control the level sets of $r_\spadesuit$ by a buckling argument, and hence prove the CLT scaling.
\end{enumerate}
In the present paper we implement this strategy for the degenerate random conductance model, which is the prototypical example in which not all algebraic moments are finite.

\subsection{Notation}
\begin{itemize} 
    \item $d\ge 1$ denotes the dimension, $|\cdot|$ denotes the Euclidean distance, $\{\vec e_i\}_{i=1,\dots,d}$ is the canonical basis of $\Z^d$;
    \item $\calz^d$ denotes the following set of (oriented) edges in $\ZZ^d$, $\calz^d:=\{ \{x,x+\vec e_i\} : x\in \ZZ^d,i=1,\dots,d\}$;
    \item For $A\subset \ZZ^d$ and $F:A\to\RR^m$ with $m\in\mathbb{N}\backslash\{0\}$, $\int_{A} F =\int_{x\in A} F(x):=\sum_{x\in A}F(x)$, $\fint_{A} F =\fint_{x\in A} F(x):=|A|^{-1}\int_{ A} F $;
    \item Similarly, for $\mathcal{A}\subset \calz^d$ and $F:\mathcal{A}\to\RR^m$  with $m\in\mathbb{N}\backslash\{0\}$, $\int_{\mathcal{A}} F=\int_{e\in \mathcal{A}} F(e):=\sum_{e\in \mathcal{A}}F(e)$, $\fint_{\mathcal{A}} F=\fint_{e\in \mathcal{A}} F(e):=|\mathcal{A}|^{-1}\int_{\mathcal{A}} F $;
    \item For $R>0$ and $x\in\ZZ^d$, $B_R(x):=\{y\in\ZZ^d\colon |x-y|\le R\}$ and $\BB_R(x):=\{e=\{y,z\} \in \calz^d:y,z\in B_R(x)\}$, we also use the shorthand notation $B_R=B_R(0)$;
    \item For $\lambda>0$ and $B=B_{r}(x)$, define $\lambda B:=B_{\lambda r}(x)$;
    \item For $f:\ZZ^d\to\RR$, $g:\calz^d\to\RR$, $ h:\ZZ^d\to\RR^d$ and $e=\{x,x+\vec e_i\}$, set $f_e:=\frac12(f(x)+f(x+\vec e_i))$, $g_e:=g(e)$ and $h_e=  h(x)\cdot\vec e_i$;
    \item For a function $f:\ZZ^d\to\RR$, the (forward) gradient $\nabla f:\ZZ^d\to\RR^d$ is defined as 
    $$\nabla f(x):=\Big(f(x+\vec{e}_i)-f(x)\Big)_{i=1,2,\cdots,d},$$
   and for $e\in\calz^d$ such that $e=\{x,x+\vec e_i\}$, $\nabla_i f(x)=\nabla_e f:=f(x+\vec e_i)-f(x)$;
    \item Similarly, the backward gradient $\nabla^* $ is defined, for $f:\ZZ^d\to\RR$, as
    $$\nabla^* f(x):=\Big(f(x)-f(x-\vec{e}_i)\Big)_{i=1,2,\cdots,d};$$
    For $f:\calz^d\to \RR$ and $ g:\ZZ^d\to\RR^d$, we set $(\nabla^*\cdot f)(x):=\sum_i \big(f(\{x,x+\vec e_i\})-f(\{x-\vec{e}_i,x\}) \big)$ and $(\nabla^*\cdot g)(x):=\sum_i\big( g(x)- g(x-\vec e_i)\big)\cdot\vec e_i$;
\item The Laplacian is defined for $f:\ZZ^d\to\RR$ by $\Delta f \coloneqq\nabla^*\cdot\nabla f$;
    \item For a $d\times d$ matrix $M$, we denote by $\|M\|$ its Frobenius norm;
    \item For any $N \ge 0$ and parameters $\alpha_1 , \ldots , \alpha_N > 0$, we write $A \lesssim_{\alpha_1 , \ldots , \alpha_N} B$ if there exists a constant $C > 0$ depending only on $\alpha_1 , \ldots , \alpha_N$ such that $A \le CB$, and $A \sim_{\alpha_1 , \cdots , \alpha_N} B$ if $A \lesssim_{\alpha_1 , \ldots , \alpha_N} B$ and $B \lesssim_{\alpha_1 , \ldots , \alpha_N} A$ -- the dependence on the dimension will be omitted for simplicity;
    \item For a real number $p\in(1,+\infty)$, $p':=\frac p {p-1}$ denotes the H\"older conjugate of $p$, that satisfies $\frac 1p + \frac 1{p'}=1$.
\end{itemize}

\section{Structure of the proof}

\subsection{Large-scale perturbative regularity}

The first object we introduce is an effective ellipticity length-scale $\rd$, that is, the random field that characterizes the minimal scale above which averages of (suitable powers of) $a$ and $a^{-1}$ are controlled.
\begin{lemma}\label{prop.ellipticr}
    There exist $C=C_d>1$ and a stationary $\frac 1{8}$-Lipschitz random field $\rd\ge 2$ such that:
    \begin{itemize}
        \item For all $x\in\ZZ^d$ and $R\ge r_\diamond(x)$,
            \begin{gather}\label{eq.equiva}
                \Big(\fint_{\BB_R(x)}a_e^{d+1}\Big)^{\frac1{d+1}} \le C\expec{a^{d+1}}^{\frac1{d+1}};\quad
\Big(\fint_{\BB_R(x)}a_e^{-(d+1)}\Big)^{\frac1{d+1}}\le C\expec{a^{-(d+1)}}^{\frac1{d+1}}.
            \end{gather}
            \item   $\rd$ satisfies the following moment bound: there exists $c_0>0$ such that
            \begin{equation}\label{eq.rdbound}
		\expec{r_\di^{c_0\gamma}}^{\frac1\gamma}\lesssim  \gamma^{\frac1{d+1}}\Gamma^{2}.
            \end{equation}
    \end{itemize}
\end{lemma}
The proof is postponed to Appendix~\ref{minimalradius}. 

Based on this, we prove quenched large-scale Meyers' estimates, from which we deduce some further elliptic regularity estimates. For simplicity, we use the shorthand notation $B_\di(x):=B_{\rd(x)}(x)$.
\begin{proposition}\label{prop.quenchedregularity}
Suppose that $u:\Z^d\to \R$ and $f:\Z^d\to \R^d$ are related on $\Z^d$ via
\begin{equation}\label{eq.u}
-\nabla^*\cdot A\nabla u=\nabla^*\cdot f.
\end{equation}
There exists $\beta>1$ depending only on $\Lambda$ such that for any $p\in[1,\beta]$ and any ball $B$, we have almost surely
\begin{equation}\label{eq.quenchedMeyers}
\fint_{x\in B}\Big(\fint_{B_{\di}(x)}A\nabla u\cdot\nabla u\Big)^{p}\,\lesssim \,\fint_{x\in 2B}\Big(\fint_{B_{\di}(x)}A^{-1} f\cdot f\Big)^{p}+\Big(\fint_{x\in 2B}\fint_{B_{\di}(x)}A\nabla u\cdot\nabla u \Big)^{p}.
\end{equation}
Furthermore, the following elliptic regularity estimates hold:
\begin{enumerate}
\item\label{cor.quenchedenergy}
(Global quenched Meyers' estimates) For $p\in[1,\beta]$,
\begin{equation}\label{eq.quenchedenergy}
\int_{x\in \ZZ^d}\Big(\fint_{B_{\di}(x)}A\nabla u\cdot\nabla u\Big)^{p}\,\lesssim  \,\int_{x\in \ZZ^d}\Big(\fint_{B_{\di}(x)}A^{-1} f\cdot f\Big)^{p}.
\end{equation}
\item\label{cor.quenchedhole-filling}
(Large-scale hole-filling estimates) For $\rd(x)\le r\le R$, if $f\equiv0$ in $B_R(x)$, then 
\begin{equation}\label{eq.quenchedhole-filling}
\int_{ B_r (x) }A\nabla u\cdot \nabla u\,\lesssim \, \Big(\frac rR\Big)^{\frac d{\beta'}}\int_{B_{R} (x) } A\nabla u\cdot \nabla u.
\end{equation}
Here we recall that $\beta':=\frac\beta{\beta-1}$ is the H\"older conjugate of $\beta$.
\end{enumerate}
\end{proposition}
With the help of these estimates, we will deduce the following annealed Meyers' estimates:
\begin{proposition}\label{prop.AnnealedM}
Suppose that $u$ and $f$ are related via~\eqref{eq.u}. There exist $\kappa, c_M>0$ depending only on $\Lambda$ such that 
\begin{enumerate}
\item\label{cor.am-ls}
For any $2\le p\le q\le2\beta$, 
\begin{equation}\label{eq.AM-ls}
\Big(\int_{x\in\ZZ^d}\expecb{\big(\fint_{B_{\di}(x)}A\nabla u\cdot \nabla u\big)^{\frac p2}}^\frac qp\Big)^{\frac1{q}} \,\lesssim_{\Lambda}\, 
\left(\int_{x\in\ZZ^d}\expecb{\Big(\fint_{B_{\di}(x)}A^{-1}f\cdot f\Big)^{\frac p2}}^\frac qp\right)^{\frac1{q}}.
\end{equation}
\item\label{cor.am}
If $\gamma \theta\ge c_M$ and $2-\kappa\le p, q, p(1+\theta)\le 2+\kappa$, then  
\begin{equation}\label{eq.AM}
\Big(\int_{\ZZ^d} \expecb{\big(A\nabla u\cdot\nabla u\big)^{\frac 12p}}^\frac qp \Big)^{\frac1q}  \lesssim_{\Lambda }  {\big(\gamma^{\frac1{d+1} } \Gamma^2\big)}^ {c_M} 
\Big(\int_{\ZZ^d}\expecb{\big(A^{-1} f\cdot f\big)^{\frac12 {p}(1+\theta)}}^\frac {q}{p(1+\theta)}\Big)^{\frac1q}.
\end{equation}
\end{enumerate}
\end{proposition}

\subsection{Spectral-gap, sensitivity calculus, and buckling}

Consider the correctors defined in Lemma~\ref{lem:liouville}. 
We now provide a first estimate on their weak norms.
\begin{proposition} \label{prop.sensitivity}
There exist $C_0,C_1,C_2\ge 1$ depending on $\Lambda$ such that, for any compactly supported $g:\ZZ^d\to\RR^{d}$ and any $1\le i,j,k\le d$, if we define $F_1(A)\coloneqq\int \nabla \phi_i \cdot g$ and $F_2(A)\coloneqq\int \nabla \sigma_{ijk} \cdot g$, then for any $\theta\in(0,\frac12)$,
\begin{equation}\label{eq.uniformbound-observables}
\expecb{|(F_1,F_2)|^{2p}}^{\frac1p}\lesssim_{ \Lambda }  (p\gamma\Gamma)^{C_0} \expecb{\Big(A (\nabla\phi_i +e_i)\cdot (\nabla\phi_i +e_i)\Big)^{p(1+\theta)}}^{\frac{1}{p(1+\theta)}} \|g\|^2_{L^2(\Z^d)},
\end{equation}
provided $\gamma\theta\ge C_1 p$ and $p\ge C_2$. 
\end{proposition}
Notice that if we take $g=\mathds1_{B_R}$ so that $\|g\|_{L^2}=R^{\frac d2}$, we recover an estimate for the constant $\Cc_R$ characterizing the CLT scaling, uniformly in $R$. This proposition can also be used to estimate the term on the right-hand side of \eqref{eq.uniformbound-observables} itself, that involves the corrector gradient. This leads us to:
\begin{proposition} \label{prop.cormoment}
There exist $C_0,C_1,C_2\ge1$ depending only on $\Lambda$, such that for $C_2\le p\le C_1^{-1}\gamma$, the (stationary) corrector gradient satisfies
\begin{equation}\label{eq.uniformbound-gradientcor}
\expecb{\Big((\nabla\phi+e)\cdot A(\nabla\phi+e)\Big)^p}^\frac1p\,\lesssim_{\Lambda}\, (p\gamma\Gamma)^{C_0}.
\end{equation}
\end{proposition}
Theorem~\ref{thm.main} follows from combining the two propositions above.

\section{Large-scale perturbative regularity}\label{sec:perturb}
In this section, we prove both the quenched and the annealed Meyers estimates of Propositions~\ref{prop.quenchedregularity} and~\ref{prop.AnnealedM}. 
Since the origin plays no special role in the argument, all balls considered below are centered at $0$ unless specified otherwise. In the first two parts of the proof, we work at the level of edges so that the discrete Caccioppoli inequality takes a neat form.

\subsection{Proof of Proposition~\ref{prop.quenchedregularity}}
The proof is divided into three parts. First, we prove the Caccioppoli-type estimate \eqref{eq.Caccioppoli} below. Second, we derive the reverse-H\"older estimate \eqref{eq.reverseholder1} and upgrade it to the quenched Meyers estimate \eqref{eq.quenchedMeyers}. Finally, we show how to obtain the remaining elliptic regularity estimates in Proposition~\ref{prop.quenchedregularity}.

In the proof, we need to manipulate double averages over balls with spatially inhomogeneous radii. For this purpose, the following lemma will be used repeatedly, and its proof is given in Appendix~\ref{apdx:averages}.
\begin{lemma}\label{lem.avrg}
    Let $B=B_r(0)$ and let $F:\ZZ^d\to\RR$ be a non-negative function.
    \begin{enumerate}[label=(\alph*),ref=\thelemma(\alph*)] 
        \item \label{lem.rgerd}
        If $r\ge2\rd(0)$, then
        \begin{equation}\label{eq.avg-equiv}
                \fint_{\frac12 B}F\,\lesssim\, \fint_{x\in B}\fint_{B_\diamond(x)}F\,\lesssim\, \fint_{2B}F.
        \end{equation}
        \item \label{lem.rlerd}
         If $r\le 2\rd(0)$, then for any $p\in(0,1)$,
        \begin{equation}\label{eq.avrg-reverseholder}
            \sup_{x\in B}\Big(\fint_{B_{{\diamond}}(x)} F\Big)^p\,\lesssim \,\fint_{x\in 5B}\Big(\fint_{B_{{\diamond}}(x)} F\Big)^p.
        \end{equation}
    \end{enumerate}
\end{lemma}
Moreover, we will use the Poincar\'e-Sobolev inequality in the following form (see  \cite[Theorem~2.6]{MourratOtto2016}).
\begin{lemma}\label{lem.sobolev}
    Let $f:\ZZ^d\to\RR$. Then for all $R\ge 1$,
    \begin{equation}\label{eq.sobolev}
        \left(\fint_{B_R}\Big|f-\fint_{B_R }f\Big|^{q}\right)^\frac1q\,\lesssim_{p}\, R\left(\fint_{ B_{R}}|\nabla f|^{p}\right)^\frac1p
    \end{equation}
    provided that $\frac1p\le\frac1q+\frac1d$.
\end{lemma}

\subsubsection{A Caccioppoli-type estimate}
We first record the standard discrete Caccioppoli inequality for solutions of \eqref{eq.u}: for any ball $B_r$ and any $c\in\RR$,
\begin{equation}\label{eq.Caccioppoli}
    \fint_{e\in\BB_r}a_e|\nabla_eu|^2\lesssim r^{-2}\fint_{e\in\BB_{ 2 r}}a_e(u- c)_e^2
    +\fint_{e\in\BB_{2 r}}a_e^{-1}f_e^2.
\end{equation}
Here, as before, $g_e:=\frac12\big(g(x)+g(y)\big)$ for $e=\{x,y\}$.

Since a solution $u$ of \eqref{eq.u} satisfies
\begin{equation}\label{eq.weaksol}
    \int_{e\in\calz^d}a_e\nabla_e u\nabla_e\vp =-\int_{e\in\calz^d} f_e\nabla_e\vp
\end{equation}
for every compactly supported test function $\vp$ on $\ZZ^d$, we choose $\vp=\eta^2(u-c)$, where $\eta$ is a cut-off for $B_r$ satisfying $0\le \eta\le 1$, $\operatorname{supp}\eta\subset B_{2r}$, $\eta\equiv 1$ on $B_r$, and $\sup_{\ZZ^d}\max(\eta,r|\nabla\eta|)\lesssim 1$.

Using the identity $2(ac-bd)=(a+b)(c-d)+(a-b)(c+d)$, we obtain the discrete Leibniz rule
$$\nabla_e(gh)=g_e\nabla_eh+h_e\nabla_eg.$$
Therefore, for any $e\in\calz^d$, Young's inequality gives
\begin{align*}
     \nabla_e \big(\eta^2 (u-c)\big)\nabla_e u\,=\,&  2(u-c)_e\eta_e\nabla_e u\nabla_e\eta+ (\eta^2)_e|\nabla_e u|^2 \\
     \,\ge\,&  \frac12 (\eta^2)_e|\nabla_e u|^2-2(u-c)_e^2|\nabla_e\eta|^2,\\
    f_e\nabla_e\big(\eta^2(u-c)\big)\,=\,&2f_e(u-c)_e\eta_e\nabla_e\eta+(\eta^2)_ef_e\nabla_eu\\
   \, \le\,& \frac14 a_e(\eta^2)_e|\nabla_eu|^2+ a_e(u-c)_e^2|\nabla_e\eta|^2+2a_e^{-1}(\eta^2)_ef_e^2.
\end{align*}
Inserting these two inequalities into \eqref{eq.weaksol} yields
\begin{equation*}
    \int_{e\in\calz^d}(\eta^2)_e a_e|\nabla_e u|^2\,\lesssim\,\int_{e\in\calz^d}(\eta^2)_e a_e^{-1}f_e^2+\int_{e\in\calz^d} |\nabla_e\eta|^2a_e(u-c)_e^2,
\end{equation*}
which is exactly \eqref{eq.Caccioppoli} by the definition of $\eta$.

\subsubsection{Quenched Meyers' estimates via Gehring's lemma}
We now prove the quenched Meyers' estimate \eqref{eq.quenchedMeyers}. To this end, we first establish a reverse-H\"older-type inequality: there exists $\frac12\le s <1$ such that
\begin{align}\label{e.Meyers-q}
    \fint_{x\in B_r} \fint_{B_{\di}(x)}A\nabla u\cdot\nabla u \,\lesssim\, \fint_{x\in B_{2r}} \fint_{B_{\di}(x)}A^{-1} f\cdot f+
\Lambda\Big(\fint_{x\in B_{2r}} \big(\fint_{B_{\di}(x)}A\nabla u\cdot\nabla u\big)^s\Big)^{\frac1s}.
\end{align}
Then we show that \eqref{eq.quenchedMeyers} follows from Gehring's lemma applied to \eqref{e.Meyers-q}.

Notice that, if $r\le 2\rd(0)$, inequality \eqref{e.Meyers-q} follows directly from \eqref{eq.avrg-reverseholder}. We therefore only consider the case $r\ge 2\rd(0)$, where \eqref{eq.Caccioppoli} is used. 

Assume $d \ge 2$. Since $r>\rd(0)$, \eqref{eq.equiva} is applicable, and thus by H\"older's inequality with indices $(d+1,1+\frac1d)$ and $(\frac{d^2+2d+2}{d},\frac{d^2+2d+2}{d^2+d+2})$ combined with Sobolev embedding (Lemma~\ref{lem.sobolev}), we have
\begin{eqnarray*}
    \inf_c r^{-2}\fint_{e\in\BB_{ 2r}}a_e\big|(u- c)_e\big|^2
&\,\le\, &r ^{-2}\Big(\fint_{e\in\BB_{ 2r}}a_e^{d+1}\Big)^{\frac 1{d+1}}\inf_c\Big(\fint_{e\in\BB_{ 2r}}\big|(u- c)_e\big|^{2+\frac2d}\Big)^{\frac {d}{d+1}}\\
   & \,\le\, & \Big(\fint_{e\in\BB_{ 2r}}a_e^{d+1}\Big)^{\frac 1{d+1}} \Big(\fint_{e\in \BB_{2r}}|\nabla_e u|^{\frac{2d(d+1)}{d^2+2d+2}}\Big)^{\frac {d^2+2d+2}{d(d+1)}}\\
    &\,\lesssim\, & \Lambda  \Big(\fint_{e\in \BB_{2r}}\big(a_e|\nabla_e u|^2\big)^{\frac{d(d+1)}{d^2+d+2}}\Big)^{\frac {d^2+d+2}{d(d+1)}}.
\end{eqnarray*}
If we set $s:=\frac{d(d+1)}{d^2+d+2}\in[\frac12,1)$, in view of the above inequality, \eqref{eq.Caccioppoli} takes the form
\begin{equation*}
    \fint_{e\in \BB_r}a_e|\nabla_e u|^2\lesssim  \fint_{e\in\BB_{2r}} a^{-1}_e|f_e|^2+\Lambda \Big(\fint_{e\in\BB_{2r}}\big(a_e|\nabla_eu|^2\big)^s\Big)^{\frac1s},
\end{equation*}
or equivalently,
\begin{equation}\label{eq.reverseholder1}
    \fint_{B_r}A\nabla u\cdot\nabla u\lesssim \fint_{B_{2r}} A^{-1} f\cdot f+\Lambda \Big(\fint_{ B_{2r}}\big(A\nabla u\cdot\nabla u\big)^s\Big)^{\frac1s}.
\end{equation}
We now quickly argue this also holds for $d=1$. By the Sobolev embedding of $W^{1,1}$ in $L^\infty$,
\begin{eqnarray*}
    \inf_c r^{-2}\fint_{e\in\BB_{ 2r}}a_e\big|(u- c)_e\big|^2
&\,\le\, &r ^{-2}\Big(\fint_{e\in\BB_{ 2r}}a_e \Big) \inf_c\sup_{e\in\BB_{ 2r}}\big|(u- c)_e\big|^{2}\\
   & \,\le\, &\Big(\fint_{e\in\BB_{ 2r}}a_e^2 \Big)^{\frac12} \Big(\fint_{e\in \BB_{2r}}|\nabla_e u|\Big)^2\\
    &\,\lesssim\, & \Lambda  \Big(\fint_{e\in \BB_{2r}}\big(a_e|\nabla_e u|^2\big)^{\frac34}\Big)^\frac {4}{3},
\end{eqnarray*}
that is, \eqref{eq.reverseholder1} with $s=\frac34$.

Hence, by the Lipschitz regularity of the field $\rd$ in the form of Lemma~\ref{lem.rgerd}, we deduce
\begin{equation*}
        \fint_{x\in B_r} \fint_{B_{\di}(x)}A\nabla u\cdot\nabla u \,\lesssim\, \fint_{x\in B_{2r}} \fint_{B_{\di}(x)}A^{-1} f\cdot f+
\Lambda \Big(\fint_{x\in B_{2r}} \fint_{B_{\di}(x)}\big(A\nabla u\cdot\nabla u\big)^s\Big)^{\frac1s}.
\end{equation*}
And \eqref{e.Meyers-q} follows as claimed by Jensen's inequality.

We then use the following version of Gehring's lemma, and a sketch of the proof is provided in Appendix~\ref{Gehring}.
\begin{lemma}\label{lem.Gehring}
    Suppose that $U,V:\ZZ^d\to\RR$ are non-negative and there exists $p>1$ such that 
    \begin{equation}\label{eq.Gehring-1}
        \fint_{B} U^p\,\le\,  \fint_{2B} V^p + C \Big(\fint_{2B}U\Big)^p 
    \end{equation}
    for a constant $C$ and all balls $B$. Then there exists $\bar p>p$ depending on $C$ and $p$ such that for any $p\le q \le \bar p$,
    \begin{equation}\label{eq.Gehring-2}
        \fint_{B} U^{q}\,\lesssim \,\fint_{2B} V^{q} + \Big(\fint_{2B}U\Big)^{q}. 
    \end{equation}
\end{lemma}
Applying this lemma to
\begin{gather*}
    U(x)=\Big(\fint_{B_{\di}(x)}A\nabla u\cdot\nabla u\Big)^{s},\quad
    V(x)=\Big(\fint_{B_{\di}(x)}A^{-1} f\cdot f\Big)^{s},
\end{gather*}
and $p=\frac1s$ implies that there exists $\tilde\beta>\frac1s$ depending only on $\Lambda$ such that \eqref{eq.Gehring-2} holds for $q\in[s^{-1},\tilde\beta]$. This is \eqref{eq.quenchedMeyers} with $\beta:=s\tilde\beta>1$.

\subsubsection{Other elliptic regularity estimates}
As a direct consequence, we can rewrite \eqref{eq.quenchedMeyers} as
\begin{equation*}
    \int_{x\in B_r}\Big(\fint_{B_{\di}(x)}A\nabla u\cdot\nabla u\Big)^{p}\lesssim  \int_{x\in B_{2r}}\Big(\fint_{B_{\di}(x)}A^{-1} f\cdot f\Big)^{p}+r^{d(1-p)}\Big(\int_{x\in B_{2r}}\fint_{B_{\di}(x)}A\nabla u\cdot\nabla u \Big)^{p}.
\end{equation*}
Sending $r\to\infty$ yields \eqref{eq.quenchedenergy} for $p\in(1,\beta]$. The inequality \eqref{eq.quenchedenergy} for $p=1$ follows from the same argument applied to \eqref{e.Meyers-q}.

Applying \eqref{eq.quenchedMeyers} to $B_R$ also provides the hole-filling estimate: if $u$ satisfies $\nabla^*\cdot A \nabla u=0$ on $B_R$, then by \eqref{eq.avg-equiv},
\begin{eqnarray*}
    \int_{ B_r}A\nabla u\cdot \nabla u\,&\lesssim&\, \int_{x\in B_{2r}}\fint_{B_\di(x)}A\nabla u\cdot \nabla u\nonumber\\
    \,&\lesssim&\, r^{\frac{d}{\beta'}} R^{\frac d\beta}\left(\fint_{x\in B_{R/4}}\Big(\fint_{B_\di(x)} A\nabla u\cdot \nabla u \Big)^{\beta}\right)^{\frac1\beta}\nonumber\\
    \,&\lesssim &\, r^{\frac{d}{\beta'}} R^{\frac d\beta} \fint_{x\in B_{R/2}} \fint_{B_\di(x)}A\nabla u\cdot \nabla u\nonumber\\
    \,&\lesssim &\, \Big(\frac rR\Big)^{\frac d{\beta'}}\int_{B_{R}} A\nabla u\cdot \nabla u,
\end{eqnarray*}
whenever $2\rd(0)\le r\le \tfrac18R$. (The statement with $r\ge \rd(0)$ then follows.)

\subsection{Proof of Proposition~\ref{prop.AnnealedM}}
The proof mainly follows the strategy of \cite{DO-20}, where the authors provide a robust method for passing from quenched estimates to annealed ones. We first show how to turn the quenched large-scale Meyers estimates \eqref{eq.quenchedenergy} into the annealed large-scale version \eqref{eq.AM-ls}. Then we argue that \eqref{eq.AM} follows from a post-processing of \eqref{eq.AM-ls} using the moment bound \eqref{eq.rdbound} of $\rd$.
\subsubsection{Annealed Meyers' estimates: large-scale version}
In this part, we prove the large-scale annealed Meyers estimates \eqref{eq.AM-ls} using the regularity estimates in Proposition~\ref{prop.quenchedregularity}.

As in \cite{DO-20}, we turn the quenched estimates into annealed estimates using the following lemma, which was first proved in \cite{CP-98}.
\begin{lemma}\label{lem.shen}
    Suppose that $1\le p_0< p_1 \le \infty$, $C,C_0\ge1$, and $g,h\in L^{p_0}\cap L^{p_1}(\ZZ^d)$ are non-negative. Assume that on any ball $B$, there exist $g_{B,0},g_{B,1}\ge0$  such that
    \begin{gather}
        g\le g_{B,0}+g_{B,1}, \quad g_{B,1}\le g+g_{B,0}\quad  \text{ on } B;\nonumber\\
        \Big(\fint_B|g_{B,0}|^{p_0}\Big)^{\frac1{p_0}} \le C\Big(\fint_{C_0B}|h|^{p_0}\Big)^{\frac1{p_0}};\label{eq.lemshen1}\\
        \Big(\fint_{\frac1{C_0}B}|g_{B,1}|^{p_1}\Big)^{\frac1{p_1}} \le C\Big(\fint_{B}|g_{B,1}|^{p_0}\Big)^{\frac1{p_0}}.\label{eq.lemshen2}
    \end{gather}
    Then for all $q\in(p_0,p_1)$,
    \begin{equation}
        \int_{\ZZ^d} |g|^q\lesssim_{p_0,p_1,q,C,C_0}\int_{\ZZ^d} |h|^q. 
    \end{equation}
\end{lemma}
To apply this lemma, given a solution $u$ of \eqref{eq.u} and a ball $B=B_r(0)$, we define $u_{B,0},u_{B,1}$ to be the solutions of the equations below:
\begin{equation}
    -\nabla^* \cdot A \nabla u_{B,0}= \nabla^*\cdot (f\mathds1_B),\quad -\nabla^* \cdot A \nabla u_{B,1}= \nabla^*\cdot (f\mathds1_{B^c}),
\end{equation}
so that $u = u_{B,0}+u_{B,1}$. In the following, we prove \eqref{eq.lemshen1} and \eqref{eq.lemshen2} with $2\le p_0=p\le p_1=2\beta$ and
\begin{gather*}
    g(x)=\expecb{\big(\fint_{B_\di(x)}A\nabla u\cdot \nabla u\big)^{\frac p2}}^\frac1p,\quad h(x)=\expecb{\big(\fint_{B_\di(x)}A^{-1}f\cdot f\big)^{\frac p2}}^\frac1p,\\
    g_{B,0}(x)=\expecb{\big(\fint_{B_\di(x)}A\nabla u_{B,0}\cdot \nabla u_{B,0}\big)^{\frac {p}2}}^{\frac 1p},\quad g_{B,1}(x)=\expecb{\big(\fint_{B_\di(x)}A\nabla u_{B,1}\cdot \nabla u_{B,1}\big)^{\frac {p}2}}^{\frac 1p}.
\end{gather*}
That is,
\begin{gather}
    \Big(\fint_{x\in  B}\expecb{\big(\fint_{ B_\di(x)}A\nabla u_{B,0}\cdot \nabla u_{B,0}\big)^{\frac {p}2}}\Big)^{\frac1{p}} \le C\Big(\fint_{x\in C_0B}\expecb{\big(\fint_{ B_\di(x)}A^{-1}f\cdot f\big)^{\frac p2}}\Big)^{\frac1{p}};\label{eq.lemshen-real1}\\
    \Big(\fint_{x\in \frac1{C_0} B}\expecb{\big(\fint_{ B_\di(x)}A\nabla u_{B,1}\cdot \nabla u_{B,1}\big)^{\frac {p}2}}^\frac{p_1}p\Big)^{\frac1{p_1}} \le C\Big(\fint_{x\in B} \expecb{\big(\fint_{ B_\di(x)}A\nabla u_{B,1}\cdot \nabla u_{B,1}\big)^{\frac {p}2}}\Big)^{\frac1{p}}.\label{eq.lemshen-real2}
\end{gather}

For \eqref{eq.lemshen-real1}, it suffices to show that
\begin{gather}
    \fint_{x\in B}\big(\fint_{ B_\di(x)}A\nabla u_{B,0}\cdot \nabla u_{B,0}\big)^{\frac {p}2} \lesssim \fint_{x\in C_0B}\big(\fint_{ B_\di(x)} A^{-1}f\cdot f\big)^{\frac p2},\label{eq.lemshen-real11}
\end{gather}
so that \eqref{eq.lemshen-real1} follows by taking expectation.
When $r>\frac14\rd(0)$, the (global) quenched Meyers estimate \eqref{eq.quenchedenergy} yields
\begin{align*}
    \fint_{x\in B}\big(\fint_{ B_\di(x)}A\nabla u_{B,0}\cdot \nabla u_{B,0}\big)^{\frac {p}2}&\,\le\,|B|^{-1}\int \big(\fint_{ B_\di(x)}A\nabla u_{B,0}\cdot \nabla u_{B,0}\big)^{\frac {p}2}\\
    &\,\lesssim \,|B|^{-1}\int \big(\fint_{ B_\di(x)} A^{-1}f\cdot f\mathds1_B\big)^{\frac {p}2}\\
    &\,\lesssim \,\fint_{x\in 6B}\big(\fint_{ B_\di(x)}A^{-1}f\cdot f\big)^{\frac {p}2}.
\end{align*}
The last inequality follows because $B\cap B_\di(x)\neq\emptyset$ and $r>\frac14\rd(0)$ imply $|x|\le r+\rd(x)\le 6r$, so that the integration region can be restricted from the whole space to $6B$.
If instead $r\le\frac14\rd(0)$, then for $x\in 2B$ one has $B\subset B_\di(x)$. Hence, by \eqref{eq.quenchedenergy},
\begin{align*}
    \fint_{ B_\di(x)} A\nabla u_{B,0}\cdot \nabla u_{B,0} &\le r_\diamond^{-d}(x)\int  A\nabla u_{B,0}\cdot \nabla u_{B,0}\\
    &\lesssim r_\diamond^{-d}(x)\int  A^{-1}f\cdot f\mathds1_{B}\lesssim  \fint_{B_\di(x)}A^{-1}f\cdot f.
\end{align*}
Therefore, \eqref{eq.lemshen-real11} follows by raising both sides to the $p/2$-th power and integrating.

For \eqref{eq.lemshen-real2}, by Minkowski's inequality and \eqref{eq.quenchedMeyers}, since $\nabla^*\cdot A\nabla u_{B,1} =0$ on $B$, one has
\begin{eqnarray*}
    &&\,\Big(\fint_{x\in \frac1{2} B}\expecb{\big(\fint_{ B_\di(x)}A\nabla u_{B,1}\cdot \nabla u_{B,1}\big)^{\frac {p}2}}^\frac{p_1}p\Big)^{\frac1{p_1}}\\
    &\,\le &\,\expecb{\Big(\fint_{x\in \frac1{2} B}\big(\fint_{ B_\di(x)}A\nabla u_{B,1}\cdot \nabla u_{B,1}\big)^{\frac {p_1}2}\Big)^{\frac p{p_1}}}^\frac{1}p\\
    &\,\lesssim &\,\expecb{\Big(\fint_{x\in B}\big(\fint_{ B_\di(x)}A\nabla u_{B,1}\cdot \nabla u_{B,1}\big)^{\frac {p}2}\Big)}^\frac{1}p,
\end{eqnarray*}
provided $2\beta\ge p_1\ge p\ge2$. Therefore, we deduce by Lemma~\ref{lem.shen} that  
$$ \Big(\int_{x\in\ZZ^d}\expecb{\big(\fint_{B_{\di}(x)}A\nabla u\cdot \nabla u\big)^{\frac p2}}^\frac qp\Big)^{\frac1{q}} \,\lesssim_{p,q,\Lambda}\, 
            \left(\int_{x\in\ZZ^d}\expecb{\Big(\fint_{B_{\di}(x)}A^{-1}f\cdot f\Big)^{\frac p2}}^\frac qp\right)^{\frac1{q}}$$
holds for $2<p\le q<2\beta$. Since the same estimate holds by \eqref{eq.quenchedenergy} for $p=q=2$, real interpolation then implies \eqref{eq.AM-ls} for all $2\le p\le q\le2\beta$, up to slightly decreasing $\beta$. 

\subsubsection{Annealed Meyers' estimates: average-free version}
To replace the large-scale averages in \eqref{eq.AM-ls} by pointwise values as in \eqref{eq.AM}, we use the probabilistic argument of \cite{DO-20}. 

The starting point is the following convex inequality: for $p,q>0$ and positive numbers $(a_i)_{i=1}^N$,
\begin{equation}\label{eq.equiv-discreteavrg}
    \Big(\frac1N\sum_{i=1}^N a_i^q\Big)^\frac1q \le N^{(\frac1p-\frac1q)_+}\Big(\frac1N\sum_{i=1}^N a_i^p \Big)^\frac1p .
\end{equation}
For a fixed $R\ge1$, $q\ge p$, and $F:\ZZ^d\to\RR$ positive, by \eqref{eq.equiv-discreteavrg} and Minkowski's inequality:
\begin{align}\label{eq.AM-l}
    \left(\int_{x\in\ZZ^d}\expecb{\Big(\fint_{ B_{R}(x)}F^2\Big)^{\frac p2}}^\frac qp\right)^{\frac1{q}}\,\ge&\, R^{-d(\frac12-\frac1q)_+}\left(\int_{x\in\ZZ^d}\expecb{\Big(\fint_{ B_{R}(x)}F^q\Big)^{\frac pq}}^\frac qp\right)^{\frac1{q}}\nonumber\\
    \,\ge&\,R^{-d(\frac12-\frac1q)_+}\left(\int_{x\in\ZZ^d}\fint_{ B_{R}(x)}\expecb{F^p}^\frac qp\right)^{\frac1{q}}\nonumber\\
    \,=&\,R^{-d(\frac12-\frac1q)_+}\left(\int_{x\in\ZZ^d}\expecb{F^p(x)}^\frac qp\right)^{\frac1{q}}.
\end{align}
Similarly, by \eqref{eq.equiv-discreteavrg} and Jensen's inequality, we obtain the following converse inequality:
\begin{align}\label{eq.AM-r}
    \left(\int_{x\in\ZZ^d}\expecb{\Big(\fint_{ B_{R}(x)}F ^2\Big)^{\frac p2}}^\frac qp\right)^{\frac1{q}}\,\le\,& R^{d(\frac1p-\frac12)_+}\left(\int_{x\in\ZZ^d}\expecb{\fint_{ B_{R}(x)}F ^p}^\frac qp\right)^{\frac1{q}}\nonumber\\
    \,\le\,&R^{d(\frac1p-\frac12)_+}\left(\int_{x\in\ZZ^d}\fint_{ B_{R}(x)}\expecb{F ^p}^\frac qp\right)^{\frac1{q}}\nonumber\\
    \,=\,&R^{d(\frac1p-\frac12)_+}\left(\int_{x\in\ZZ^d}\expecb{F ^p(x)}^\frac qp\right)^{\frac1{q}}.
\end{align}
Now we want to replace $R$ by $\rd$ in the two inequalities above, and consider the case where $q\ge p\ge 2$. By Jensen's inequality applied to the function $x\mapsto x^{\frac qp}$, we have
\begin{align*}
    \int_{x\in\ZZ^d}\expecb{F^p(x) }^\frac qp
    \,\lesssim \,& \int_{x\in\ZZ^d}\expecb{\sum_{k=1}^{+\infty}2^{-kd }\rd(x)^{d }\mathds1_{2^{k-1}\le\rd(x)<2^{k}}F ^{p}(x)}^\frac qp\\
    \,\lesssim \,& \sum_{k=1}^{+\infty}2^{kd (\frac qp-1)} \int_{x\in\ZZ^d}\expecb{\mathds1_{2^{k-1}\le\rd(x)<2^{k}}F ^{p}(x)}^\frac qp.
\end{align*}
We then apply  \eqref{eq.AM-l} to the integral on the right-hand side and obtain
\begin{align*}
    &\int_{x\in\ZZ^d}\expecb{F^p(x) }^\frac qp
    \,\lesssim \, \sum_{k=1}^{+\infty}2^{kd (\frac qp-1)} 2^{kdq(\frac12-\frac1q)_+}\int_{x\in\ZZ^d}\expecb{\Big(\fint_{B_{2^{k-2}}(x)}\mathds1_{2^{k-1}\le\rd<2^{k}}F^2\Big) ^{\frac p2}}^\frac qp.
\end{align*}
Since $\rd$ is $\frac18$-Lipschitz, $B_{2^{k-2}}(x)\cap\{2^{k-1}\le \rd<2^k\}\neq \emptyset$ implies $2^{k-2}\le \rd(x)<2^{k+1}$. Thus we further deduce that
\begin{align*}
    \int_{x\in\ZZ^d}\expecb{F^p(x) }^\frac qp
    \,\lesssim \,& \sum_{k=1}^{+\infty}2^{kd (\frac qp-1)} 2^{kdq(\frac12-\frac1q)_+}\int_{x\in\ZZ^d}\expecb{\mathds1_{2^{k-2}\le\rd(x)<2^{k+1}}\Big(\fint_{B_\diamond(x)} F^2\Big) ^{\frac p2}}^\frac qp\\
    \, {\lesssim} \,& \int_{x\in\ZZ^d}\sum_{k=1}^{+\infty}\expecb{ r_\diamond^{dp\big( (\frac 1p-\frac1q)+(\frac12-\frac1q)_+\big)}   \mathds1_{2^{k-2}\le\rd(x)<2^{k+1}}\Big(\fint_{B_{\di}(x)}F^2\Big) ^{\frac p2}}^\frac qp.
\end{align*}
Since $\ell^1\hookrightarrow \ell^\frac qp$, the above inequality implies, by H\"older's inequality with indices $(\frac rp,\frac{r}{r-p})$ for $r>p$,
\begin{align*}
    \int_{x\in\ZZ^d}\expecb{F^p(x) }^\frac qp
    \,\lesssim \,&  \int_{x\in\ZZ^d} \expecb{ r_\diamond^{dp\big( (\frac 1p-\frac1q)+(\frac12-\frac1q)_+\big)}    \Big(\fint_{B_{\di}(x)}F^2\Big) ^{\frac p2}}^\frac qp\\
    \,\lesssim \,& \expec{r_\diamond^{d\frac{pr}{r-p}\big( (\frac 1p-\frac1q)+(\frac12-\frac1q)_+\big)}}^{q(\frac1p-\frac1r)} \int_{x\in\ZZ^d} \expecb{   \Big(\fint_{B_{\di}(x)}F^2\Big) ^{\frac r2}}^\frac qr.
\end{align*}
In the last inequality, we have used the stationarity of $\rd$. Hence, since $2\le p\le q$, provided
\begin{equation}\label{eq.consld1}
   2d\frac{pr}{r-p} (\frac12-\frac1q )\le   c_0    \gamma ,
\end{equation}
by the moment bound \eqref{eq.rdbound} of $\rd$, we obtain 
\begin{equation}\label{eq.functional1}
    \Big(\int_{x\in\ZZ^d}\expecb{F^p(x) }^\frac qp\Big)^\frac1q \, \lesssim \, {\big(\gamma^{\frac1{d+1} } \Gamma^2\big)}^{\frac{2d}{c_0} (\frac12-\frac1q )} \Big(\int_{x\in\ZZ^d} \expecb{   \Big(\fint_{B_{\di}(x)}F^2\Big) ^{\frac r2}}^\frac qr\Big)^\frac1q.
\end{equation}

\noindent
Similarly, by \eqref{eq.AM-r} , for $r>p\ge2$,
\begin{align*}
     \int_{x\in\ZZ^d}\expecb{\Big(\fint_{B_{\di}(x)}F^2\Big) ^{\frac {p}2}}^\frac q{p} 
    \,&\lesssim \,  \sum_{k=1}^{+\infty}2^{kd (\frac qp-1)} 2^{kdq(\frac1p-\frac12)_+}\int_{x\in\ZZ^d}\expecb{ \mathds1_{2^{k-1}\le\rd<2^{k}}F^p}^\frac qp\\
    \,&\lesssim \,   \expec{r_\diamond^{d\frac{pr}{r-p} (\frac 1p-\frac1q) }}^{q(\frac1p-\frac1r)} \int_{x\in\ZZ^d} \expecb{  F (x)^ r}^\frac qr.
\end{align*}
Hence, again, by  the moment bound \eqref{eq.rdbound}, we obtain
\begin{equation}\label{eq.functional2}
     \Big(\int_{x\in\ZZ^d}\expecb{\Big(\fint_{B_{\rd}(x)}F^2\Big) ^{\frac {p}2}}^\frac q{p}\Big)^\frac1q \,\lesssim\,    {\big(\gamma^{\frac1{d+1} } \Gamma^2\big)}^{ \frac d{c_0}(\frac1p-\frac1q) } \Big(\int_{x\in\ZZ^d}\expecb{F^{r}(x)}^\frac q{r}\Big)^\frac1q,
\end{equation}
provided
\begin{equation}\label{eq.consld2}
    d\frac{pr}{r-p} (\frac 1p-\frac1q) \le   c_0  \gamma  .
\end{equation}

Applying the two functional inequalities \eqref{eq.functional1} and \eqref{eq.functional2} to \eqref{eq.AM-ls} yields
\begin{equation}\label{eq.AM-q>p}
    \Big(\int\expecb{\big(A\nabla u\cdot \nabla u\big)^{\frac p{2(1+\delta)}}}^\frac {(1+\delta)q}p\Big)^{\frac1{q}} \lesssim_{\Lambda}  {\big(\gamma^{\frac1{d+1} } \Gamma^2\big)}^ {\frac{3d}{c_0}|\frac12-\frac1q|} 
    \Big(\int\expecb{\big(A^{-1}f\cdot f\big)^{\frac p{2(1-\delta)}}}^\frac {(1-\delta)q}p\Big)^{\frac1{q}}.
\end{equation} 
for $2\le p\le q \le 2\beta$ provided $\gamma\delta\ge  \frac{2d}{c_0}p(\frac12-\frac1q)$. A duality argument yields the above inequality for $2\ge p\ge q\ge(2\beta)'$ and $\gamma\delta\ge \frac{2d}{c_0}p(\frac1q-\frac12)$. Real interpolation
entails that the above inequality holds for $|\frac12-\frac1p|\lor|\frac12-\frac1q|\le \frac13(\frac12-\frac1{2\beta})$ and $\gamma\delta\ge \frac{2d}{c_0}( \beta - 1 )$, with the multiplicative constant replaced by $\big({\gamma^{\frac1{d+1}}} \Gamma^2\big)^ {\frac{3d}{2c_0}(1-\frac1{\beta})}$. 
Finally, set $\kappa:=\frac12(1-\frac1\beta)$, by  choosing suitable $p,q,\delta$,  \eqref{eq.AM} follows under the hypothesis that $2-\kappa\le p,q,p(1+\theta)\le 2+\kappa$ and that  $\gamma\theta\ge \frac{4d}{c_0}(\beta-1)=:c_M$.

\section{Proof of Proposition~\ref{prop.sensitivity} and Proposition~\ref{prop.cormoment}}
In this section, we prove Proposition~\ref{prop.sensitivity} and Proposition~\ref{prop.cormoment} by first considering the case where $a$ is replaced by a truncated version $a_M:=(a\land M)\lor M^{-1}$ ($M\ge 1$), so that every expectation appearing here is a priori finite thanks to \cite{GO1}. Since the estimates obtained in this section are independent of $M$, we shall argue by approximation and take the limit $M\uparrow+\infty$.

In the following, we use the spectral-gap inequality to quantify ergodicity; a proof can be found in \cite[Lemma~2.3]{GO1} for i.i.d. conductances.
\begin{lemma}
    For any measurable random variable $X=X(A)$,
    \begin{equation}\label{eq.SG-2}
        \textup{Var}(X)\le \frac12 \expecb{ \int_{x\in \ZZ^d} |D_x X |^2},
    \end{equation}
    where 
    \begin{equation*}
     D_x X(A):=X(A)-X(A^{(x)}).
    \end{equation*}
Here $A^{(x)}$ is the random field defined by $A^{(x)}(x')=A(x')$ for $x'\neq x$, while $A^{(x)}(x)$ is an i.i.d. copy of $A(x)$.
    Moreover, this implies, for $p\ge 1$,
    \begin{equation}\label{eq.SG-p}
        \expecb{|X-\expec X|^{2p}}^{\frac1p}\lesssim 4 p^2\expecb{\Big(\int_{x\in \ZZ^d} |D_x X |^2\Big)^p}^\frac1p.
    \end{equation}
\end{lemma}

This section is divided into three parts. In the first part, we prove Proposition~\ref{prop.sensitivity} by deducing a handy representation formula for (functional) derivatives of  $(F_1,F_2)$ and then showing that they can be controlled by the corrector gradient. In the second part, we prove Proposition~\ref{prop.cormoment} using Proposition~\ref{prop.sensitivity} and a buckling argument. Finally, we present the approximation argument that removes the uniform ellipticity assumption in the limit $M\uparrow +\infty$.

\subsection{Proof of Proposition~\ref{prop.sensitivity}}
In this section, we first derive the representation formulas for $F_1=F_1(A)$ and $F_2=F_2(A)$, which depend on $A$ through \eqref{eq.corrector} and \eqref{eq.flux}. The dependence of $(F_1,F_2)$ on the directions $i,j,k$ is left implicit.
\subsubsection{Representation formulas for  functional derivatives }
In this part, we prove the following identity:
\begin{gather}
    D_x F_1=\Big((A-A^{(x)})(\nabla\phi_i ^{(x)}+e_i)\cdot\nabla u \Big)(x);\label{eq.representationformula1}\\
    D_x F_2=\Big((A-A^{(x)})(\nabla\phi_i^{(x)}+e_i)\cdot (\nabla^*_j v e_k-\nabla^*_k v e_j+\nabla w_{jk})\Big)(x), \label{eq.representationformula2}
\end{gather}
where $u$, $v$, $w$ are given by
\begin{gather}
	-\nabla^*\cdot A \nabla u=\nabla^*\cdot g,\label{eq.ug}\\
	-\Delta v = \nabla^*  \cdot g,\label{eq.v}\\
	-\nabla^*\cdot A  \nabla w_{jk}=  \nabla ^* \cdot A  (\nabla^*_k v e_j-\nabla^*_j v e_k)\label{eq.w}.
\end{gather}

In fact, given $A$ and $A'$ two realizations of the random field,
\begin{equation*}
    F_1(A)-F_1(A')=\int_{\Z^d}(A-A')\big(\nabla\phi_i (A')+e_i\big) \cdot \nabla u.
\end{equation*}
This is because, by definition,
\begin{gather*}
    -\nabla^* \cdot A\nabla\phi_i (A)=\nabla^*\cdot Ae_i; \\
    -\nabla^* \cdot A'\nabla\phi_i (A')=\nabla^*\cdot A'e_i;\\
\end{gather*}
one has 
\begin{equation}\label{eq.error}
-\nabla^*\cdot A\nabla[\phi_i (A)-\phi _i(A')]=\nabla^*\cdot(A-A')[\nabla\phi _i(A')+e_i].
\end{equation}
By the weak formulation of this equation,
\begin{align*}
    F_1(A)-F_1(A')&=\int_{\Z^d}\nabla[\phi_i (A)-\phi_i (A')] \cdot g\\
    &=-\int_{\Z^d} A\nabla[\phi_i (A)-\phi_i(A')] \cdot  \nabla u\\
    &=\int_{\Z^d} (A-A')(\nabla\phi_i (A')+e_i)\cdot \nabla u.
\end{align*}
Thus, by taking $A'=A^{(x)}$, we obtain \eqref{eq.representationformula1}. For \eqref{eq.representationformula2}, set $q_i(A)\coloneqq A(\nabla\phi_i(A)+e_i)$. By definition \eqref{eq.flux} of $\sigma$, one has
\begin{eqnarray*}
  F_2(A)-F_2(A')& = &\int_{\Z^d} \big (\nabla\sigma_{ijk}(A)-\nabla\sigma_{ijk}(A')\big)\cdot  g\\
    & = & -  \int_{\Z^d} \big(\nabla\sigma_{ijk}(A)-\nabla\sigma_{ijk}(A')\big)\cdot  \nabla v \\
  & = &   \int_{\Z^d} (e_j\otimes e_k-e_k\otimes e_j)\big( q_i(A)-q_i(A') \big)\cdot  \nabla^* v.
\end{eqnarray*}
Since 
\begin{equation*}
  q_i(A)-q_i(A') =  A\big(\nabla \phi_i(A)-\nabla\phi_i(A')\big)+ (A-A')(\nabla\phi_i(A')+e_i),
\end{equation*}
one has, by \eqref{eq.w} and \eqref{eq.error},
\begin{eqnarray*}
& &  F_2(A)-F_2(A')\\
& = & - \int A\big(\nabla \phi_i(A)-\nabla\phi_i(A')\big)\cdot  \nabla w_{jk} + \int (A-A')(\nabla\phi_i(A')+e_i)\cdot (\nabla^*_j v e_k-\nabla^*_k v e_j)\\
& = &  \int (A-A')(\nabla\phi_i(A')+e_i)\cdot (\nabla^*_j v e_k-\nabla^*_k v e_j+\nabla w_{jk}).
\end{eqnarray*}
Therefore, \eqref{eq.representationformula2} follows by taking $A'=A^{(x)}$.

\subsubsection{Control of $F$ by corrector gradient}

In this subsection, we prove that
\begin{equation}\label{eq.cont-r1}
    \expecb{(F_1,F_2)^{2p}}^{\frac1p}\lesssim_{\Lambda}p^2 {\gamma^{\frac{ 2c_M}{d+1}}}\Gamma^{ 4c_M  +3} \expecb{\Big(A(\nabla\phi_i+e_i)\cdot (\nabla\phi_i+e_i)\Big)^{p(1+\theta)}}^{\frac1{p(1+\theta)}}  \|g\|_{L^2(\Z^d)}^2,
\end{equation}
provided that $\theta\gamma\ge C_1 p$ and $p\ge  C_2$ for constants $C_1,C_2\ge1$ depending on $\Lambda$.

By \eqref{eq.SG-p} applied to the centered random variable $F$ and \eqref{eq.representationformula1}--\eqref{eq.w}, one has
\begin{align*}
    \expecb{F^{2p}}^\frac1p\,\lesssim&\,   p^2 \expecb{\Big(\int_{x\in\ZZ^d}\big((A-A^{(x)})(\nabla\phi_i^{(x)}+e_i) \cdot \nabla u\big)^2\Big)^p}^\frac1p\\
&+p^2\expecb{\Big(\int_{x\in\ZZ^d}\big((A-A^{(x)})(\nabla\phi_i^{(x)}+e_i) \cdot(\nabla^*_j v e_k-\nabla^*_k v e_j)\big)^2\Big)^p}^\frac1p\\
&+p^2\expecb{\Big(\int_{x\in\ZZ^d}\big((A-A^{(x)})(\nabla\phi_i^{(x)}+e_i) \cdot \nabla w_{jk}\big)^2\Big)^p}^\frac1p.
\end{align*}
We only treat the term containing $w$, since it requires solving two equations at a time. The other two terms can be controlled in a similar way. By $L^p$-$L^{p'}$ duality, one has
\begin{align*}
&\expecb{\Big(\int_{x\in\ZZ^d}\big((A-A^{(x)})(\nabla\phi_i^{(x)}+e_i) \cdot \nabla w_{jk}\big)^2\Big)^p}^\frac1p\\
    =& \sup_{Y\colon\expec{|Y|^{2p'}}=1}\expecb{\int_{x\in\ZZ^d}\big((A-A^{(x)})(\nabla\phi_i^{(x)}+e_i) \cdot \nabla (Yw_{jk})\big)^2}.
\end{align*}
Since 
\begin{align*}
    &\Big((A-A^{(x)})(\nabla\phi_i^{(x)}+e_i) \cdot \nabla (Y w_{jk})\Big)^2 \\
    \le & \|A^{-\frac12}(A-A^{(x)})(A^{(x)})^{-\frac12}\|^2\Big(A^{(x)}(\nabla\phi_i^{(x)}+e_i)\cdot (\nabla\phi_i^{(x)}+e_i)\Big)\Big(A\nabla (Yw_{jk})\cdot\nabla (Yw_{jk})\Big),
\end{align*}
H\"older's inequality with indices $\big( \frac{p'}{\tilde\theta}, p(1+\theta) , \frac{p'}{1+\tilde\theta} \big)$ (with $\theta,\tilde\theta>0$ constants to be fixed later) together with the fact that $(A^{(x)},\phi^{(x)}_i)$ and $(A,\phi_i)$ are identically distributed yields that 
\begin{align*}
    &\expecb{ \int_{x\in\ZZ^d}\big((A-A^{(x)})(\nabla\phi_i^{(x)}+e_i) \cdot \nabla (Yw_{jk})\big)^2}\\
    \, \le \,&\expecb{\|A^{-\frac12}(A-A^{(\cdot)})(A^{(\cdot)})^{-\frac12}\|^\frac{2p'}{\tilde\theta}}^{\frac{\tilde\theta}{p'}}\expecb{\Big(A(\nabla\phi_i+e_i)\cdot (\nabla\phi_i+e_i)\Big)^{p(1+\theta)}}^{\frac{1}{p(1+\theta)}}\\
&\phantom{\expecb{\|A^{-\frac12}(A-A^{(x)})(A^{(x)})^{-\frac12}\|^\frac{2p'}{\tilde\theta}}^{\frac{\tilde\theta}{p'}}\expec\expec} \cdot \int_{\ZZ^d} \expecb{\Big(A\nabla (Yw_{jk})\cdot\nabla (Yw_{jk})\Big)^\frac{p'}{1+\tilde\theta}}^\frac{1+\tilde\theta}{p'}.
\end{align*}
Notice that in this case, $\theta$ and $\tilde\theta$ are related via $2\tilde\theta(p-1)=\frac\theta{1+\theta}$. The first term on the right-hand side is bounded by $\Gamma^2$ provided $\frac{2p'}{\tilde\theta}\le\gamma $. By the annealed Meyers estimates~\eqref{eq.AM} applied to \eqref{eq.w}, we have
\begin{align*}
&\int_{\ZZ^d} \expecb{\Big(A\nabla (Yw_{jk})\cdot\nabla (Yw_{jk})\Big)^\frac{p'}{1+\tilde\theta}}^\frac{1+\tilde\theta}{p'}\\
\,\lesssim_{\Lambda } \,&  (\gamma^{\frac{1}{d+1}}\Gamma^{ 2 })^{2c_M} \int_{\ZZ^d}\expecb{\Big( A\nabla^* (Yv)\cdot\nabla^* (Yv)\Big)^\frac{p'}{1+\tilde\theta/2}}^\frac{1+\tilde\theta/2}{p'}.
\end{align*}
Therefore, H\"older's inequality with indices $(p',\frac{2p'}{\tilde\theta})$, combined with the energy estimate for the deterministic equation~\eqref{eq.v}, yields
\begin{align*}
\int_{\ZZ^d}\expecb{\Big( A\nabla^* (Yv)\cdot\nabla^* (Yv)\Big)^\frac{p'}{1+\tilde\theta/2}}^\frac{1+\tilde\theta/2}{p'}\,&\le\, \expecb{\Big(\|A\|Y^2\Big)^\frac{p'}{1+\tilde\theta/2}}^\frac{1+\tilde\theta/2}{p'}\int_{\ZZ^d}|\nabla^* v|^2\\
\,&\le\, \expecb{\|A\|^\frac{2p'}{\tilde\theta}}^\frac{\tilde\theta}{2p'}\int_{\ZZ^d}|g|^2.
\end{align*}
Hence \eqref{eq.cont-r1} follows, provided that $\frac{2p'}{\tilde\theta}\le\gamma $ and that \eqref{eq.AM} is applicable. The first constraint reads
$$\frac{2p'}{\tilde\theta}\le\gamma \quad \Leftarrow\quad 8p \le \gamma\theta,$$
and the application of the annealed Meyers estimates requires that 
\begin{equation}\label{eq.cons-appMeyers}
    1+\tilde\theta\le p'\le 1+\kappa/2, \quad\frac{2\gamma\tilde\theta}{2+\tilde\theta}\ge c_M.
\end{equation}
These conditions translate into $p\ge 1+2\kappa^{-1}\eqqcolon C_2$ and $C_1 p\le \gamma\theta$ for $C_1\coloneqq 4c_M \lor 8$.

\subsection{Proof of Proposition~\ref{prop.cormoment}}
To this end, define the stationary random field
\begin{equation*}
r_\spd(x):=\inf\Big\{r\ge\rd(x) \textup{ dyadic }:  \forall R\ge r \textup{ dyadic}, \forall i, \fint_{e\in\BB_R(x)}a_e|\nabla_e\phi_i|^2\le C_\spd\fint_{e\in\BB_{2R}(x)}a_e \Big\}.
\end{equation*}
To establish estimates on $r_\spd$, we further define 
\begin{equation*}
r_{\spd,i}(x):=\inf\Big\{r\ge\rd(x) \textup{ dyadic }:  \forall R\ge r \textup{ dyadic},  \fint_{e\in\BB_R(x)}a_e|\nabla_e\phi_i|^2\le C_\spd\fint_{e\in\BB_{2R}(x)}a_e \Big\},
\end{equation*}
and note that $r_{\spd} = \max_{1\le i \le d}r_{\spd,i}$.
In what follows, we implicitly fix $i$ and skip the subscript $i$ of $r_{\spd,i}$.
The random radius $r_\spd(x)$ is almost surely finite by the ergodic theorem if we have $C_\spd > \expec{a |\nabla \phi_i|^2}/\expec{a }$.
Notice that \eqref{eq.uniformbound-gradientcor} follows from the moment bound on $r_\spd$ because, by the definition of $r_\spd$ and \eqref{eq.quenchedhole-filling},
\begin{align}\label{eq.cont-ed-rspd}
\Big(A(\nabla\phi_i+e_i)\cdot (\nabla\phi_i+e_i)\Big)(0)\,\le&\, \int_{ B_{\di}} A(\nabla\phi_i+e_i)\cdot (\nabla\phi_i+e_i) \nonumber\\
\,\lesssim &\, \Big(\frac{\rd}{r_\spd}\Big)^{\frac d{\beta'}}r_\spd^d\fint_{ B_{r_\spd}}A(\nabla\phi_i+e_i)\cdot (\nabla\phi_i+e_i)\nonumber\\
\,\lesssim &\, \Gamma r_\di^{\frac d{\beta'}}r_\spd^{\frac d{\beta}}\le \Gamma r_\spd^d.
\end{align}
\noindent
Hence, the $dp$-th moment of $ r_\spd$ controls the $p$-th moment of $A(\nabla\phi_i+e_i)\cdot (\nabla\phi_i+e_i)$. Provided we control moments of $r_\spd$, we obtain the proposition by relabelling the exponents.
To estimate $r_\spd$, it suffices to control its level sets. Our starting point is the following set inclusion: for dyadic $R\ge2$,
\begin{equation}\label{eq.setdecomp}
\{r_\spd=R\} \subset \{R/2 < \rd <R\} \cup \{r_\spd=R \ge 2\rd\}.
\end{equation}
As the probability of the first set on the right-hand side is already controlled by \eqref{eq.rdbound}, we only treat the second set in the following.

To this end, we first notice that, by definition,  the  inequalities
\begin{gather}
    \fint_{e\in\BB_R}a_e|\nabla_e\phi_i|^2\le C_\spd\fint_{e\in\BB_{2R}}a_e;\label{eq.spddef1}\\
    \fint_{e\in\BB_{R/2}}a_e|\nabla_e\phi_i|^2\ge C_\spd\fint_{e\in\BB_{R}}a_e
    \label{eq.spddef2}
\end{gather}
hold on  $\{r_\spd=R>2\rd\}$. 
Combined with Caccioppoli's inequality~\eqref{eq.Caccioppoli} (whose constant is denoted $C$ here), \eqref{eq.spddef1} yields that, for any $c\in\RR$,
\begin{align*}
    C_\spd\fint_{e\in\BB_{R}}a_e\,\le\,\fint_{e\in\BB_{R/2}}a_e|\nabla_e\phi_i|^2\,\le\, C R^{-2}\fint_{e\in\BB_{R}}a_e(\phi_i-c)_e^2+C\fint_{e\in\BB_{R}}a_e.
\end{align*}
Therefore, if $C_\spd\ge  2C$, on  $\{r_\spd=R>2\rd\}$
\begin{equation*}
    \fint_{e\in\BB_{R}}a_e\lesssim  R^{-2}\fint_{e\in \BB_{R}}a_e(\phi_i-c)_e^2.
\end{equation*}
And H\"older's inequality with indices $(d+1,1+\frac1d)$ combined with Jensen's inequality implies that
\begin{align*}
    1\,\lesssim &\,R^{-2}\Big(\fint_{e\in\BB_R}a_e\Big)^{-1} \Big(\fint_{e\in\BB_R}a_e^{d+1}\Big)^{\frac1{d+1}} \Big(\fint_{e\in\BB_R}(\phi_i-c)_e^{2(1+\frac 1d )}\Big)^{\frac d{d+1}}\\
    \,\lesssim&\,
R^{-2} \Lambda\Big(\fint_{B_R}(\phi_i-c)^{2(1+\frac 1d )}\Big)^{\frac d{d+1}}.
\end{align*}
We then use the following lemma that first appeared in \cite{BellaKniely}. A proof is given in Appendix~\ref{sobolev}.
\begin{lemma}\label{lem.avrgsobolev}
    Let $S,s\ge 1$, $\tau:=d\big(\frac1s-\frac1S\big)\in[0,1]$ and $\mu\in (0,1)$. Then for any function $\psi:\ZZ^d\to \RR$,
    \begin{equation}\label{eq.sob-avrg}
        R^{-1}\Big(\fint_{B_R}|\psi-\fint_{B_R}\psi|^S\Big)^{\frac1S}\lesssim_s R^{-(1-\tau)(1-\mu)}\Big(\fint_{B_{2R}}|\nabla\psi|^s\Big)^{\frac1s}+\Big(\fint_{B_R}\Big|\fint_{B_{R^\mu}(x)}\nabla\psi\Big|^s\Big)^{\frac1s}.
    \end{equation}
\end{lemma}
The application of this lemma with $S=2(1+\frac1d)$, $s=\frac{2(d+1)}{d+2}$, $\tau=\frac d{d+1}\in(0,1)$ and $\mu$ to be fixed later yields that, on $\{r_\spd=R>2\rd\}$,
\begin{align*}
    \Lambda^{-\frac12}\,&\lesssim\,  R^{-\frac{1-\mu}{d+1}}\Big(\fint_{B_{2R}}|\nabla \phi_i|^s\Big)^{\frac1s}+\Big(\fint_{x\in B_R}\Big|\fint_{B_{R^\mu}(x)}\nabla\phi_i\Big|^s\Big)^{\frac1s}\nonumber\\
    \,&\lesssim \,R^{-\frac{1-\mu}{d+1}}\Big(\fint_{e\in\BB_{2R}}a_e|\nabla_e \phi_i|^2\Big)^{\frac12}\Big(\fint_{e\in\BB_{2R}}a_e^{-(d+1)}\Big)^{\frac1{2(d+1)}}+\Big(\fint_{x\in B_R}\Big|\fint_{B_{R^\mu}(x)}\nabla\phi_i\Big|^s\Big)^{\frac1s}\nonumber\\
    \,&\lesssim \Lambda^{\frac12}\, R^{-\frac{1-\mu}{d+1}}+\Big(\fint_{x\in B_R}\Big|\fint_{B_{R^\mu}(x)}\nabla\phi_i\Big|^2\Big)^{\frac12}.
\end{align*}
In the last two steps we have used H\"older's inequality with exponents $(\frac {d+2}{d+1},d+2)$, Jensen's inequality, \eqref{eq.equiva}, and the definition of $r_\spd$ at scale $2R$. Therefore, if $R\ge R_0:=R_0(\mu,\Lambda)$ is chosen large enough, the first term on the right-hand side can be absorbed into the left-hand side. Thus, on $\{r_\spd=R>2\rd\}$, for any $p>1$,
\begin{equation*}
    1\,\lesssim_{ \Lambda}\,   \fint_{x\in B_R}\Big|\fint_{B_{R^\mu}(x)}\nabla\phi_i\Big|^2\,\lesssim_{ \Lambda}\,  \Big( \fint_{x\in B_R}\Big|\fint_{B_{R^\mu}(x)}\nabla\phi_i\Big|^{2p}\Big)^\frac1p.
\end{equation*}
 Hence, since the law of $\nabla\phi_i$ is invariant under translation, we may take $p$-th moment on both sides of the above inequality to get 
\begin{equation*}
    \expecb{\mathds1_{r_\spd=R>2\rd}}^\frac1p \,\lesssim_{ \Lambda}\,  \expecb{ \Big( \fint_{x\in B_R}\Big|\fint_{B_{R^\mu}(x)}\nabla\phi_i\Big|^{2p}\Big)}^\frac1p=\expecb{  \Big|\fint_{B_{R^\mu} }\nabla\phi_i\Big|^{2p} }^\frac1p.
\end{equation*}
Applying \eqref{eq.cont-r1} to $g=|B_{R^\mu}|^{-1}\mathds1_{B_{R^\mu}}$ and using \eqref{eq.cont-ed-rspd} yields for $\gamma\theta \ge  (8c_M\lor 16) p$ and $p\ge 1+\frac2\kappa$,
\begin{align*}
    \expecb{\mathds1_{r_\spd=R>2\rd}}^\frac1p  &\,\lesssim_\Lambda\,     p^2 {\gamma^{ \frac {2c_M} {d+1} }}\Gamma^{  4c_M +3 }
\expecb{\Big(A(\nabla\phi_i+e_i)\cdot (\nabla\phi_i+e_i)\Big)^{p(1+\theta/2)}}^{\frac1{p(1+\theta/2)}}  R^{-d\mu }\\
    &\,\lesssim_\Lambda\, p^2 {\gamma^{ \frac {2c_M} {d+1} } }\Gamma^{  4c_M +4 }R^{-d\mu }\expecb{\Big(r_\di^{\frac d{\beta'}}r_\spd^{\frac d{\beta}}\Big)^{p(1+\theta/2)}}^{\frac1{p(1+\theta/2)}}.
\end{align*}
By H\"older's inequality with indices $\big(\frac{1+\theta}{\theta/2},\frac{1+\theta}{1+\theta/2}\big)$, this entails
\begin{align}\label{eq.bucklingbase}
    \expecb{\mathds1_{r_\spd=R>2\rd}}^\frac1p   &\,\lesssim_\Lambda\, p^2 {\gamma^{ \frac {2c_M} {d+1} }}\Gamma^{  4c_M +4 }R^{-d\mu }\expecb{r_\di^{\frac {pd}{\beta'}(1+\theta)(1+2/\theta)}}^{\frac1{p(1+\theta)(1+2/\theta)}} \expecb{r_\spd^{  {\frac {pd}{\beta}}(1+\theta)}}^{\frac1{p(1+\theta)}}\nonumber\\
    &\,\lesssim_\Lambda\, p^2 {\gamma^{ \frac {2c_M+ d/(c_0\beta')} {d+1} }}\Gamma^{  4c_M +4 +2d/(c_0\beta')}R^{-d\mu }  \expecb{r_\spd^{  {\frac {pd}{\beta}}(1+\theta)}}^{\frac1{p(1+\theta)}},
\end{align}
and we need that 
\begin{equation}\label{cons.rd}
    \frac{pd}{\beta'}(1+\theta)(1+2/\theta)\le c_0 \gamma\quad\Leftarrow\quad \theta\gamma\ge 2c_Mp
\end{equation} 
in the last inequality in order to control $\rd$ by \eqref{eq.rdbound}. Here we have used that by definition, $c_M= \frac{4d}{c_0}(\beta-1)\ge 4d/(c_0\beta')$. Thus, one has
\begin{equation}\label{eq.bucklingbase1}
\expecb{\mathds1_{r_\spd=R>2\rd}} \le C_\Lambda^p p^{2p} {\gamma^{ \frac {3c_M p } {d+1} }}\Gamma^{  (5c_M +4)p}R^{-d\mu p}  \expecb{r_\spd^{  {\frac {pd}{\beta}}(1+\theta)}}^{\frac1{(1+\theta)}},
\end{equation}
where $C_\Lambda$ is the multiplicative constant in \eqref{eq.bucklingbase}.

Finally, applying the layer-cake formula and inserting \eqref{eq.bucklingbase1} into \eqref{eq.setdecomp} entails
\begin{align*}
    \expecb{r_\spd^{pd(1-\frac1{2\beta'})}}&\,\le\, R_0^{pd(1-\frac1{2\beta'})}+\sum_{m\in\N:2^m\ge R_0} 2^{mpd(1-\frac1{2\beta'})}\expecb{ \mathds1 _{r_\spd=2^m} }\\
&\,\le\, R_0^{pd(1-\frac1{2\beta'})}+\sum_{m\in\N:2^m\ge R_0}2^{mpd(1-\frac1{2\beta'})}  \Big( \expecb{ \mathds1 _{2^{m-1}<\rd\le2^m} }  + \expecb{ \mathds1 _{r_\spd=2^m>2\rd} }\Big)\\
&\,\lesssim\, R_0^{pd(1-\frac1{2\beta'})}+ \expecb{r_\di^{pd(1-\frac1{2\beta'})}} \\
&\phantom{\,\lesssim\, R_0 } +C_\Lambda^p p^{2p} {\gamma^{ \frac {3c_M p } {d+1} }}\Gamma^{  (5c_M +4)p}\expecb{r_\spd^{  {\frac {pd}{\beta}}(1+\theta)}}^{\frac1{p(1+\theta)}} \sum_{m\in\N:2^m\ge R_0}2^{mpd(1-\frac1{2\beta'}-\mu)}   .
\end{align*}
We then take $\mu=1-\frac1{4\beta'}$, so that $1-\frac1{2\beta'}-\mu=-\frac1{4\beta'}<0$. In this case, the sum converges and is bounded by a constant.  Furthermore, let $\theta=\frac{1}{2\beta'}$, which ensures that $\beta^{-1}(1+\theta)<1-\frac1{2\beta'}$. This implies that the expectation $\expecb{r_\spd^{  {\frac {pd}{\beta}}(1+\theta)}}^{\frac1{p(1+\theta)}}$ in the last term on the right-hand side can be absorbed into the left-hand side by Young's inequality:
\begin{align*}
    \expecb{r_\spd^{pd(1-\frac1{2\beta'})}}^{\frac1{pd(1-\frac1{2\beta'})} }&\,\lesssim_{\Lambda} \,  \expecb{r_\di^{pd(1-\frac1{2\beta'})}}^{\frac1{pd(1-\frac1{2\beta'})} } +  p^{ \frac {2\beta'}d } {\gamma^{ \frac {4c_M \beta' } {d^2  } }}  \Gamma^{  \frac {2\beta'}d (5c_M +4) }.
\end{align*}
The first term is bounded by \eqref{eq.rdbound} if $pd(1-\frac1{2\beta'})\le c_0\gamma$. In the end, all the constraints translate into 
\begin{equation}\label{cond.p}
    C_2\coloneqq 1+\frac2\kappa\le p \le C_1^{-1}  \gamma,
\end{equation}
where $C_1\coloneqq (16c_M\beta')\lor (32\beta')\lor (d/c_0)$. And there exists another constant $C_0>0$  depending only on $\Lambda$ such that  $\expecb{r_\spd^{p }}^{\frac 1{p}}\,\lesssim_\Lambda \,   (p\gamma\Gamma)^{C_0}$, provided that  $p$ satisfies \eqref{cond.p}. 

\subsection{Approximation argument}
Set $\phi_{M}:=\phi(A_M)$. We shall prove that $\expec{|\nabla\phi-\nabla\phi_M|^r}\to 0$ for $1\le r\le\frac{2\gamma}{\gamma+1}$, so that the (uniform) estimates \eqref{eq.uniformbound-observables} and \eqref{eq.uniformbound-gradientcor} that yield moment bounds for $r_\spd$ are retained when $M$ tends to infinity. In fact, for any $M\ge1$
\begin{align*}
    \expec{|\nabla\phi -\nabla\phi_M|^r}\le&\, \expec{|\nabla\phi -\nabla\phi_M|^{\frac{2\gamma}{\gamma+1}}}^{\frac{\gamma+1}{2\gamma}r}\\
\le&\, \expec{A (\nabla\phi -\nabla\phi_M)\cdot (\nabla\phi -\nabla\phi_M)}^{\frac r2}\expec{\|A^{-1}\| ^{\gamma}}^{\frac {r}{2\gamma}}.
\end{align*}
By~\eqref{eq.error} as well as~\eqref{eq.cont-ed-rspd}, we can control the right-hand side:
\begin{align*}
    &\expec{A (\nabla\phi -\nabla\phi_M)\cdot (\nabla\phi -\nabla\phi_M)}\\
\,\le\,&\sum_{i=1}^d\expecb{\|A ^{-\frac12}(A -A_{M})A_{M}^{-\frac12} \|^2\big(A_{M }(\nabla\phi_{M }+e_i)\cdot(\nabla\phi_{M }+e_i)\big)}\\
    \,\lesssim\,&\expecb{\|A ^{-\frac12}(A -A_{M})A_{M}^{-\frac12} \|^2  r_\spd^d(A_{M})}\\
	\,\le\,&\expecb{\|A ^{-\frac12}(A -A_{M})A_{M}^{-\frac12} \|^\gamma}^{\frac2\gamma}\expecb{r_\spd^{\frac {d\gamma}{\gamma-2}}(A_{M})}^{1-\frac2\gamma}.
\end{align*}
Since $\frac{|a-b|}{\sqrt{ab}}\le a +a^{-1}+b +b^{-1}$, we see that the right-hand side tends to $0$  by dominated convergence theorem when $M\to+\infty$, provided $\expecb{r_\spd^{\frac {d\gamma}{\gamma-2}}(A_M)}^{1-\frac2\gamma}$ is uniformly bounded in $M$.
Therefore, as long as $\gamma\ge 2+ C_1 d$, the condition is verified, and $\expec{|\nabla\phi-\nabla\phi_M|^r}\to 0$.

We now give the corresponding approximation argument for $\sigma$.
Recall that $\delta_M\sigma:=\sigma_{ijk}(A)-\sigma_{ijk}(A_M)$ satisfies, with $\delta_M \phi:=\phi_i(A)-\phi_i(A_M)$ and $\delta_M A:=A-A_M$,
\begin{equation*}
-\Delta \delta_M\sigma \,=\, \nabla_{j} \big(\delta_MA  (\nabla\phi_{i}(A)+{e_i})+A_M \nabla\delta_M\phi\big)_k-\nabla_{k} \big(\delta_MA  (\nabla\phi_{i}(A)+{e_i})+A_M \nabla\delta_M\phi\big)_j.
\end{equation*}
Let $1<r= \frac{2\gamma}{\gamma+1}<2$ (which is possible since $\gamma>1$). By Calder\'on-Zygmund theory, the equation yields
\begin{equation*}
\expec{|\nabla \delta_M\sigma|^r}^\frac1r \lesssim \expec{|\delta_MA  (\nabla\phi_{i}(A)+{e_i})|^r}^\frac1r +\expec{|A_M \nabla\delta_M\phi|^r}^\frac1r.
\end{equation*}
By H\"older's inequality with exponents $(\gamma+1,\frac{\gamma+1}\gamma)$, we obtain
\[
 \expec{|\delta_MA  (\nabla\phi_{i}(A)+{e_i})|^r}^\frac1r \le  \expecb{\|A^{-\frac12}(\delta_MA)^2 A^{-\frac12}\|^\gamma}^\frac1{2\gamma} \expec{A(\nabla\phi_{i}(A)+{e_i})\cdot (\nabla\phi_{i}(A)+{e_i})}^\frac12,
 \]
and
\[
\expec{|A_M \nabla\delta_M\phi|^r}^\frac1r\le  \expecb{\|A^{-\frac12}A_M^2 A^{-\frac12}\|^\gamma}^\frac1{2\gamma} \expec{A \nabla \delta_M \phi \cdot \nabla \delta_M \phi}^\frac12.
\]
As above, we then conclude that $\expec{|\nabla \delta_M\sigma|^\frac{2\gamma}{\gamma+1} }\to 0$.

\section{Proof of Corollary~\ref{cor.stationary} }
In this section we prove Corollary~\ref{cor.stationary}. For simplicity, we set $p\coloneqq C_1^{-1}\gamma$.
In the case $d=1$, without loss of generality, suppose that $x>0$. By \eqref{eq.uniformbound-observables}, with $g_x(y):=\mathds1_{0\le y<x}$,
\begin{equation*}
    \expec{|\phi(x)-\phi(0)|^p}^\frac1p=\expecb{\Big|\int_{\ZZ} \nabla\phi\cdot g_x\Big|^p}^\frac1p\lesssim_\Lambda (\gamma\Gamma)^{C_0} \|g_x\|_{L^2}=(\gamma\Gamma)^{C_0} \sqrt x.
\end{equation*}
Since the corrector gradient is explicitly given by $\nabla_e\phi=\expec{a_e^{-1}}^{-1} a_e^{-1}-1$, the Central Limit Theorem confirms that this bound on the spatial growth of $\phi$ is sharp. We also notice that in this case the corrector $\sigma\equiv 0$.

Now we turn to the case $d\ge 2$. In the following, we only prove the statement for $\phi$, since the estimate for $\sigma$ is obtained in the same way. Let $G$ denote the lattice Green function defined by 
\begin{equation*}
    -\Delta_yG(y,x)=\delta_x.
\end{equation*}
(Since we only need differences of Green functions, the case of dimension $d=2$ is fine using that, unlike $G$ itself, $\nabla G$ is uniquely defined.)
When $d\ge2$, it satisfies $|\nabla^2G(0,x)|\lesssim (1+|x|)^{-d}$, which implies an estimate for later use:
\begin{equation}\label{est.gradientG}
|\nabla\big(G(y , x)-G(y , 0)\big)|\lesssim \left\{
\begin{aligned}
 &\frac{|x|}{(1+|y|)^{d}},&\quad& |y|\ge 2|x|;\\
 &\frac{1}{(1+|x-y|)^{d-1}}+ \frac{1}{(1+|y|)^{d-1}}, &\quad& |y|\le 2|x|.
\end{aligned}\right.
\end{equation}
By definition, we formally have
\begin{align}\label{eq.representation-phi}
    \phi(x)-\phi(0)=\int_{\ZZ^d}\nabla\phi\cdot\nabla\big(G(\cdot , x)-G(\cdot , 0)\big).
\end{align}
This should be understood as 
\begin{equation*}
    \phi(x)-\phi(0)=\lim_{k\to\infty}\int_{\ZZ^d}\eta_k\nabla\phi\cdot\nabla\big(G(\cdot , x)-G(\cdot , 0)\big),
\end{equation*}
with $\eta_k$ a cut-off function of $B_{2^k}$, that is, $0\le\eta_k\le1$, $\eta_k\equiv1$ on $B_{2^k}$, $\operatorname{supp}\eta_k\subset B_{2^{k+1}}$, $\sup_{\ZZ^d}\max(\eta_k,2^k|\nabla\eta_k|)\lesssim1$.
This is because, by integration by parts,
\begin{align*}
    &\int_{\ZZ^d}\eta_k\nabla\phi\cdot\nabla\big(G(\cdot , x)-G(\cdot , 0)\big)\\
    \,=\,&-\int_{\ZZ^d}\phi\nabla^*\cdot\Big(\eta_k\nabla\big(G(\cdot , x)-G(\cdot , 0)\big)\Big)\\
    \,=\,&\phi(x)-\phi(0)-\sum_{\calz^d\ni e=\{y,y+\vec e_i\}}\phi(y+\vec e_i)\nabla_e\eta_k\nabla_e\big(G(\cdot , x)-G(\cdot , 0)\big).
\end{align*} 
By definition of $\eta_k$, $\nabla\eta_k$ is supported on $B_{2^{k+1}}\backslash B_{2^k}$, so that in the above sum, $|y| \ge 2^k$. When $k$ is so large that $|y|\ge 2|x|$, \eqref{est.gradientG} together with the bound $|\nabla\eta_k|\lesssim 2^{-k}$ yields
\begin{align*}
    \sum_{\calz^d\ni e=\{y,y+\vec e_i\}}\phi(y+\vec e_i)\nabla_e\eta_k\nabla_e\big(G(y , x)-G(y, 0)\big)
    \,\lesssim\,2^{-k}|x|\fint_{B_{2^{k+1}}}|\phi|.
\end{align*}
The right-hand side tends to $0$ by the sub-linearity of the corrector \eqref{est.sublinear} when $k\to+\infty$.
 
Therefore, by  Theorem~\ref{thm.main},
\begin{align*}
    \expec{|\phi(x)-\phi(0)|^p}^\frac1p\,=&\,\expecb{\Big|\int\nabla\phi \cdot\nabla\big(G(\cdot , x)-G(\cdot , 0)\big)\Big|^p}^\frac1p\\
    \,\lesssim_\Lambda\,& (\gamma\Gamma)^{C_0} \,\|\nabla G(\cdot , x)-\nabla G(\cdot , 0)\|_{L^2(\Z^d)}.
\end{align*}
The right-hand side is uniformly bounded in $x$ when $d\ge3$ since the $L^2$ norm of $y\mapsto |\nabla G|(y,x)\,\lesssim\,(1+|x-y|)^{1-d}$ is finite and independent of $x$. In the case $d=2$, by \eqref{est.gradientG}, we have $\|\nabla G(\cdot,x)-\nabla G(\cdot,0)\|_{L^2}\,\lesssim\, \sqrt{\log(1+|x|)}$, which is \eqref{eq.growth-corrector}.

When $d\ge 3$, we have  $\sup_{x}\expec{|\phi(x)|^p}^\frac1p\,\lesssim_{\Lambda,\gamma,\Gamma}\, 1$. The existence of the stationary corrector will follow from this uniform bound via a standard compactness argument in the probability space using massive correctors, see the proof in \cite{GNO-reg}.

\appendix

\section{Proof of auxiliary results}

\subsection{Control of the effective ellipticity length-scale: Proof of Lemma~\ref{prop.ellipticr}}\label{minimalradius}
Define
    \begin{equation*}
        \tilde r_\diamond^{\pm}(x):=\inf\{r \text{ dyadic } \ge 2: \forall R> r \textup{ dyadic}, \Big|\fint_{e\in\BB_R(x)}a_e^{\pm (
d+1)}-\expec{a^{\pm (d+1)}}\Big|\le \frac12\expec{a^{\pm (d+1)}}\}.
    \end{equation*}
We only prove the moment bound for $\tilde r_\diamond^{+}$, since the argument for $\tilde r_\diamond^{-}$ is identical. Let $X_R:=\fint_{e\in\BB_R}a_e^{d+1}$. By the spectral-gap inequality \eqref{eq.SG-p}, when $p(d+1)\le \gamma$,
\begin{align*}
    \expecb{\big|X_R-\expec{X_R}\big|^p}&\le  2^p p^{ p}\expecb{\Big|\int_{e\in\BB_R} C^2R^{-2d}(a_e^{d+1}-(a_e^{(e)})^{d+1})^2\Big|^{\frac p2}}\\
    &\lesssim  2^p p^{ p} R^{-\frac12pd} \expecb{a^{p(d+1)}}.
\end{align*}
Therefore, by the definition of $\tilde r_\di^+$ and Chebyshev's inequality with $p(d+1)=\gamma$,
\begin{eqnarray*}
    \mathbb P(\tilde r_\diamond^+= 2^m)\,&\le&\, \mathbb P\left(\Big|\fint_{e\in\BB_{2^m}} a_e^{d+1}-\expec{a^{d+1}} \Big|\ge \frac12\expec{a^{d+1}} \right)\\
    \,&\lesssim &\, (4p)^p \expec{a^{d+1}}^{-p} \expecb{\big|X_{2^{m}}-\expec{X_{2^{m}}} \big|^p}\\
    \,&\lesssim &\, (4 p)^p \expec{a^{d+1}}^{-p} \expec{a^{p(d+1)}}  2^{-\frac12{m}pd}\\
   \, &\lesssim & \,  (4  p)^p\Gamma^{2p(d+1)}  2^{-\frac12mpd}.
\end{eqnarray*}
The choice $p(d+1)=\gamma$ entails that for $c_0\coloneqq\frac{d}{3(d+1)}$, one has $c_0 \gamma -\frac d2p=- \frac{d}{6(d+1)} \gamma<0$. Therefore, the layer-cake formula gives
\begin{eqnarray*}
     \expec{\left(\tilde r_\diamond^+\right)^{c_0 \gamma }}\,&=&\,  \sum_{m=1}^\infty\,   2^{c_0m\gamma  }\mathbb P(\tilde r_\diamond^+= 2^m)\\
\,&\le&\,  \sum_{m=1}^\infty\,  (4 p)^p\Gamma^{2p(d+1)}  2^{c_0m\gamma -\frac d2pm}\\
   \, &\lesssim & \,  \Big( {\frac{4\gamma}{d+1}}\Big)^{\frac\gamma{d+1}}\Gamma^{2\gamma}.
\end{eqnarray*}
Hence
\begin{equation*}
	\expec{\left(\tilde r_\diamond^+\right)^{ c_0 \gamma }}^\frac1\gamma\,\lesssim \, \gamma^{\frac1{d+1}}\Gamma^{2}.
\end{equation*}

Define
\begin{equation*}
    r_\diamond(x):= \inf_y\{(\tilde r_\diamond^{+}\lor\tilde r_\diamond^{-})(y)+\frac1{8}|x-y|\}.
\end{equation*}
Then $r_\diamond$ is the maximal $\frac1{8}$-Lipschitz stationary random field bounded above by $\tilde r_\diamond^{+}\lor\tilde r_\diamond^{-}$, and hence satisfies the same bound as \eqref{eq.rdbound}. It remains to prove \eqref{eq.equiva}.

Indeed, for any $\rho\ge \rd(x)$, there exists $y$ such that $\tilde r^+_\diamond(y) + \frac18|x-y|\le \rho$. Hence $\BB_\rho(x)\subset \BB_{9\rho}(y)$ and, if $k\in\mathbb{N}$ is chosen so that $2^{k-1}\le 9\rho\le2^k$, then
\begin{align*}
    \fint_{e\in\BB_\rho(x)} a_e^{d+1}\le 9^d\fint_{e\in\BB_{9\rho}(y)} a_e ^{d+1}\le 18^d\fint_{e\in\BB_{2^k}(y)} a_e^{d+1}\le 2\cdot18^d  \expec{a^{d+1}}.
\end{align*}
Similarly,
\begin{equation*}
    \fint_{e\in\BB_\rho(x)} a_e^{-(d+1)}\le 2\cdot18^d  \expec{a^{-(d+1)}},
\end{equation*}
which completes the proof of \eqref{eq.equiva}.

\subsection{Playing with local averages: proof of Lemma~\ref{lem.avrg}}\label{apdx:averages}
In this section, we prove Lemma~\ref{lem.avrg}, following the arguments in \cite[Lemma~6.5]{DO-20}.  

\subsubsection{Proof of Lemma~\ref{lem.rgerd}}
The starting point of the proof is the observation that if $|x-y|\le \rd(x)$, then $|\rd(x)-\rd(y)|\le \frac18\rd(x)$, hence $\frac89\rd(y)\le\rd(x)\le\frac87\rd(y)$. 

The first inequality in \eqref{eq.avg-equiv} follows from the fact that 
$$\{x,y\colon y\in  \frac12 B,x\in \frac89 B_\diamond(y)\}\subset\{x,y\colon x\in B, y\in B_\diamond(x)\}.$$
This set inclusion holds because $|x-y|\le \frac89 r_\diamond(y)$ implies $y\in B_\di(x)$ and if moreover $|y|\le \frac12 r$, then $|x|\le|x-y|+|y|\le \frac89\rd(0) +  \frac{10}9 |y|\le r$. Therefore, by Fubini's theorem,
\begin{align*}
    \fint_{x\in B}\fint_{B_{\diamond}(x)}F\,=&\,r^{-d}\int_{x\in\ZZ^d}\int_{y\in\ZZ^d}\rd(x)^{-d}F(y)\mathds1_{x\in B,y\in B_\diamond(x)}\\
    \,\gtrsim&\,r^{-d}\int_{x\in\ZZ^d}\int_{y\in\ZZ^d}\rd(y)^{-d}F(y)\mathds1_{x\in B,y\in B_\diamond(x)}\\
    \,\gtrsim&\,r^{-d}\int_{x\in\ZZ^d}\int_{y\in\ZZ^d}\rd(y)^{-d}F(y)\mathds1_{y\in  \frac12 B,x\in \frac89B_\diamond(y)}\\
    \,\gtrsim&\,\fint_{\frac12 B}F.
\end{align*}
By a similar calculation, the second inequality in \eqref{eq.avg-equiv} follows from
$$\{ x\in B,y\in B_\di(x)\}\subset \{ y\in 2B,x\in\frac87 B_\di(y)\}.$$
The above holds for the reason that $|x|\le r,|x-y|\le\rd(x)$ imply $x\in\frac87B_\diamond(y)$ and $|y|\le |x|+|x-y|\le r+\frac18r+\rd(0)\le 2r$.

\subsubsection{Proof of Lemma~\ref{lem.rlerd}}
In the case $r\gtrsim \rd(0)$, the result is trivial, and we may thus assume $\rd(0)\gg r$ in the rest of the proof, in order to ensure that
$\rd(0)-\frac{15}8r>0$. 
We first observe that since $x\in B$, $B_{{\diamond}}(x)\subset {B_{\rd(0)+\frac98 r}(0)}$. This is because for all $y\in B_{{\diamond}}(x)$, $|y|\le |y-x|+|x|\le \rd(x)+r \le \rd(0)+\frac98 r$. Therefore, since $r\le 2\rd(0)$,
\begin{align*}
    \fint_{B_{{\diamond}}(x)} F \lesssim \fint_{B_{\rd(0)+\frac98 r}(0)} F.
\end{align*}
We therefore pick $N$ points $(z_i)_{i=1}^N\subset 4B$ with $N$ depending only on $d$ such that $4B\subset \cup_{i=1}^N B_{r}(z_i)$. In this way, $B_{\rd(0)+\frac98 r}(0)\subset \cup_{i=1}^N B_{\rd(0)- \frac{15}{8} r}(z_i)$. Indeed, for $y\in B_{\rd(0)+\frac98 r}(0)$, we set $\tilde y:=4r\frac{y}{|y|}\in 4B$, so that $\inf_i|y-z_i|\le|y-\tilde y|+\inf_i |\tilde y-z_i|$. By definition, $\inf_i |\tilde y-z_i|\le r$. Hence, 
since $\rd(0)\gg r$,
$$\inf_i|y-z_i|\le ||y|-4 r|+r\le \rd(0)- \frac{15}8 r.$$ 
Therefore, by the subadditivity of $t\mapsto t^p$,
\begin{equation*}
    \Big(\fint_{B_{{\diamond}}(x)} F \Big)^p \,\lesssim\, \Big(\sum_{i=1}^N\fint_{B_{\rd(0)-\frac{15}8r}(z_i)} F\Big)^p    \,\le\, \sum_{i=1}^N\Big(\fint_{B_{\rd(0)-\frac{15}8r}(z_i)} F\Big)^p.
\end{equation*}
Finally, we note that for any $y\in B_{r}(z_i)$, $B_{\rd(0)-\frac{15}8r}(z_i)\subset B_{{\diamond}}(y)$. Indeed, for $z\in B_{\rd(0)-\frac{15}8r}(z_i)$, $|z-y|\le|z-z_i|+|z_i-y|\le \rd(0)-\frac{15}8r+r\le  \rd(y)$, because $|\rd(0)-\rd(y)|\le\frac18|y|\le\frac58r$. Hence, by $\cup_{i=1}^N B_{r}(z_i)\subset 5B$,
\begin{align*}
    \Big(\fint_{B_{{\diamond}}(x)} F \Big)^p \,\lesssim&\, \sum_{i=1}^N\Big(\fint_{B_{\rd(0)-\frac{15}8r}(z_i)} F\Big)^p\\
    \,\lesssim&\, \sum_{i=1}^N\fint_{y\in B_{r}(z_i)}\Big(\fint_{B_{{\diamond}} (y)} F\Big)^p\,\lesssim \,\fint_{y\in5B}\Big(\fint_{B_{{\diamond}}(y)} F\Big)^p,
\end{align*}
which is \eqref{eq.avrg-reverseholder}.

\subsection{Gehring's Lemma: proof of Lemma~\ref{lem.Gehring}}\label{Gehring}
We proceed as in \cite{Iwaniec1998}. First, we show that \eqref{eq.Gehring-1} implies the global estimate $\|U\|_{L^q}\,\lesssim\,\|V\|_{L^q}$ for $q$ slightly larger than $p$. Then we explain how to localize this estimate in order to obtain \eqref{eq.Gehring-2}.  

\subsubsection{Global estimates with improved integrability}
For a non-negative function $g:\ZZ^d\to \RR$, define its ($p$-)maximal function by 
\begin{equation*}
    \maxi_p g(x):=\sup_{B\ni x} \Big(\fint_B g^p\Big)^\frac1p, \quad \maxi g:=\maxi_1 g.
\end{equation*}
If we assume $U,V \in L^{1}(\ZZ^d)$, taking the supremum on both sides of \eqref{eq.Gehring-1} yields $\maxi_p U\le \maxi_p V+C\maxi U $.

By classical arguments based on Vitali's covering lemma and the Calder\'on-Zygmund decomposition (which are still valid on $\ZZ^d$), one has for $g\in L^1(\ZZ^d)$ and $t>0$,
\begin{equation*}
    t^{-1}\int_{g\ge \frac12 t} g\sim |\{\maxi g\ge t\}|.
\end{equation*}
Multiplying both sides by $q t^{q-1}$ (for $q>1$) and integrating over $(0,+\infty)$ yields $2^q\|g\|_{L^q}^q\sim (1-\frac1q)\|\maxi g\|_{L^q}^q$, as long as one of them is finite. Equivalently, this reads, for $q>p$, $ \|g\|_{L^q} \sim (1-\frac pq)^{\frac1q}\|\maxi_p g\|_{L^q} $. Hence, for $q>p$,
\begin{align*}
   \| U\|_{L^q}&\lesssim  (1-\frac pq)^{\frac 1q} \|\maxi_p U\|_{L^q} \\
    &\le (1-\frac pq)^{\frac 1q}\big( \|\maxi_p V\|_{L^q}+ C\|\maxi U\|_{L^q}\big)\\
    &\lesssim   \| V\|_{L^q}+ C \Big(\frac{q-p}{q-1}\Big)^{\frac1q}\| U\|_{L^q} 
\end{align*}
Denote by $\tilde C$ the multiplicative constant in the above estimate. When $q-p$ is so small that $C\tilde C  \Big(\frac{q-p}{q-1}\Big)^{\frac1q}\le \frac12$, the second term on the right-hand side can be absorbed into the left-hand side, and we obtain $\|U\|_{L^q}\lesssim  \|V\|_{L^q}$ .

\subsubsection{Local estimates with increasing support}
To obtain a local version of the above inequality, we define, for a fixed ball $B_0$,
\begin{equation*}
    U_{B_0}(x)=\rho_{B_0}(x)U(x),\quad V_{B_0}(x)=\rho_{B_0}(x)V(x)+\Big(\int_{2B_0} U \Big)\mathds1_{2B_0}
\end{equation*}
with $\rho_{B_0}(x):=\operatorname{dist}(x,(2B_0)^c)$. Using \eqref{eq.Gehring-1} together with a geometric argument, one can show that for any ball $B$,
\begin{equation*}
    \Big(\fint_{ B} |U_{B_0}|^{p}\Big)^{\frac1p}\lesssim\Big(\fint_{2 B} |V_{B_0}|^p\Big)^{\frac1p}+C\fint_{2 B} U_{B_0}.
\end{equation*}
Then by the previous argument (notice that $U_{B_0}$ and $V_{B_0}$ are compactly supported), we obtain the global estimate $\|U_{B_0}\|_{L^q}\lesssim\| V_{B_0}\|_{L^q}$ provided $q-p\ll1$. Removing the cutoffs $\rho_{B_0}$, we obtain that 
\begin{equation*}
    \Big(\fint_{B}U^{q}\Big)^{\frac1q}\lesssim\Big(\fint_{2B}V^q\Big)^{\frac1q}+\fint_{2B}U.
\end{equation*}

\subsection{A variant of Poincar\'e-Sobolev inequality: Proof of Lemma~\ref{lem.avrgsobolev}}\label{sobolev}
In this section, we use the Poincar\'e-Sobolev inequality \eqref{eq.sobolev} to prove Lemma~\ref{lem.avrgsobolev} following the proof of \cite[Lemma~2.3]{BellaKniely}:

We introduce the intermediate scale $R^\mu$, then apply \eqref{eq.sobolev} to $B_{R^\mu}(x)$ and $B_R$ to get
\begin{align*}
    &\inf_{c}R^{-1}\Big(\fint_{B_R}|\psi-c|^S\Big)^{\frac1S}\\
    \lesssim\,&R^{-1}\Big(\fint_{x\in B_R}\fint_{B_{R^\mu}(x)}\Big|\psi-\fint_{B_{R^\mu}(x)}\psi\Big|^S\Big)^{\frac1S}+\inf_{c}R^{-1}\Big(\fint_{x\in B_R}\Big|\fint_{B_{R^\mu}(x)}\psi-c\Big|^S\Big)^{\frac1S}\\
    \lesssim_s\,& R^{\mu-1}\Big(\fint_{x\in B_R}\big(\fint_{B_{R^\mu}(x)}|\nabla\psi|^s\big)^{\frac Ss}\Big)^{\frac1S} + \Big(\fint_{x\in B_{R}}\Big|\fint_{B_{R^\mu}(x)}\nabla\psi\Big|^s\Big)^{\frac1s}.
\end{align*}
Notice that applying the Poincar\'e-Sobolev inequality \eqref{eq.sobolev} requires $\frac1s-\frac1S\le\frac1d$, or equivalently, $\tau\le1$. We then bound the first term on the right-hand side using Young's convolution inequality:
\begin{align*}
    \Big(\fint_{x\in B_R}\big(\fint_{ B_{R^\mu}(x)}|\nabla\psi|^s\big)^{\frac Ss}\Big)^{\frac1S}
    =\,&\Big(\fint_{x\in B_R}\big(\fint_{ B_{R^\mu}(x)}|{\bf 1}_{B_{2R}}\nabla\psi|^s\big)^{\frac Ss}\Big)^{\frac1S}\\
    \le\,&R^{-\frac dS-\mu \frac ds}\Big(\int_{\ZZ^d}\big({\bf 1}_{B_{R^\mu}(\cdot)}*|{\bf 1}_{ B_{2R}}\nabla\psi|^s\big)^{\frac Ss}\Big)^{\frac1S}\\
    \le\,&R^{-\frac dS-\mu \frac ds}R^{\mu  \frac dS}\Big(\int_{\ZZ^d}|{\bf 1}_{B_{2R}}\nabla\psi|^s\Big)^{\frac1s}\\
    \le\,& R^{\tau(1-\mu)}\Big(\fint_{B_{2R}}|\nabla\psi|^s\Big)^{\frac1s},
\end{align*}
which requires $S\ge s\Leftrightarrow \tau\ge0$. Hence, \eqref{eq.sob-avrg} follows.

\section*{Acknowledgements}
The authors acknowledge financial support from the European Research Council (ERC) under the European Union's Horizon 2020 research and innovation programme (Grant Agreement n$^\circ$~864066).

\bibliographystyle{abbrv}
\def\cprime{$'$}

\end{document}